\documentclass[aos]{imsart}

\RequirePackage{amsthm,amsmath,amsfonts,amssymb,newclude}
\RequirePackage[numbers,sort&compress]{natbib}
\RequirePackage[colorlinks,citecolor=blue,urlcolor=blue]{hyperref}
\RequirePackage{graphicx,xcolor}

\startlocaldefs
\theoremstyle{plain}

\newtheorem{theorem}{Theorem}[section]
\newtheorem{proposition}[theorem]{Proposition}
\newtheorem{lemma}[theorem]{Lemma}
\newtheorem{remark}[theorem]{Remark}
\theoremstyle{remark}
\newtheorem{definition}[theorem]{Definition}

\newcommand{\R}{\mathbb{R}}

\newcommand{\Vvert}{\vert \hspace{-1pt} \vert \hspace{-1pt} \vert}

\endlocaldefs

\usepackage{xr}

\begin{document}
	
	\begin{frontmatter}
		\title{PCA for Point Processes}
		\runtitle{PCA for Point Processes}
		
				\begin{aug}

						\author[A]{\fnms{Franck}~\snm{Picard}\ead[label=e1]{franck.picard@ens-lyon.fr}\orcid{0000-0001-8084-5481}}
						\author[B]{\fnms{Vincent}~\snm{Rivoirard}\ead[label=e2]{Vincent.Rivoirard@dauphine.fr}\orcid{0000-0001-6461-806X}}
											\author[B]{\fnms{Angelina}~\snm{Roche}\ead[label=e3]{roche@ceremade.dauphine.fr}\orcid{0000-0001-8490-4995}}	\and
						\author[C]{\fnms{Victor}~\snm{Panaretos}\ead[label=e4]{victor.panaretos@epfl.ch}\orcid{0000-0002-2442-9907}}
						\address[A]{Laboratoire de Biologie et Mod\'elisation de la Cellule, CNRS ENS de Lyon, \printead[presep={,\ }]{e1}}
						\address[B]{CEREMADE, Universit\'e Paris Dauphine, \printead[presep={,\ }]{e2,e3}}
						\address[C]{Institut de Math\'ematiques, Ecole Polytechnique F\'ed\'erale de Lausanne, \printead[presep={,\ }]{e4}}
					\end{aug}

		\begin{abstract}
		We introduce a novel statistical framework for the analysis of replicated point processes that allows for the study of point pattern variability at a population level. By treating point process realizations as random measures, we adopt a functional analysis perspective and propose a form of functional Principal Component Analysis (fPCA) for point processes. The originality of our method is to base our analysis on the cumulative mass functions of the random measures which gives us a direct and interpretable analysis. Key theoretical contributions include establishing a Karhunen-Lo\`{e}ve expansion for the random measures and a Mercer Theorem for covariance measures. We establish convergence in a strong sense, and introduce the concept of principal measures, which can be seen as latent processes governing the dynamics of the observed point patterns. We propose an easy-to-implement estimation strategy of eigenelements for which parametric rates are achieved. We fully characterize the solutions of our  approach to Poisson and Hawkes processes  and validate our methodology via simulations and diverse applications in seismology, single-cell biology and neurosiences, demonstrating its versatility and effectiveness. Our method is implemented in the \texttt{pppca} R-package.
		\end{abstract}
		
		\begin{keyword}[class=MSC]
			\kwd[Primary ]{62M}
			\kwd{62P}
			\kwd[; secondary ]{60G55}\kwd{62H25}
		\end{keyword}
		
		\begin{keyword}
			\kwd{Point processes}
			\kwd{Karhunen-Lo\`eve expansion}
			\kwd{PCA}
		\end{keyword}
		
	\end{frontmatter}
\section{Introduction}

Point processes constitute a ubiquitous framework that is essential in probability and statistics as well as in numerous application fields. In the most general sense, point processes are discrete random sets in an arbitrary space, serving as the natural mathematical formalism for discrete random patterns. Depending on the ambient space, point processes can model spatial/temporal events, tessellations of space, or random geometrical configurations. Poisson processes are arguably fundamental \cite{Kingman}, due to their simplicity, and because they serve as building blocks for more elaborate point process models. Key models include renewal, marked, cluster, or doubly stochastic point processes \cite{DV1,DV2}, and geometrical aspects have also long been reviewed \cite{chiu2013stochastic}.  When the ambient space is the real line, point processes represent sequences of discrete events in time and are referred to as temporal point processes, among which Hawkes processes \cite{Hawkes71a, Hawkes71b} are crucial for enhancing the flexibility of modeling temporal events, introducing a tractable and interpretable model of temporal dependency \cite{BM1, BM2}. The Poisson and Hawkes processes are among the most fundamental, popular, and analytically/statistically tractable, serving as basic examples, including in the present work. Given their formidable flexibility in handling particular characteristics of the data, the scope of application of point process models is vast, ranging from neurosciences \cite{MR961117, JMM, BONNET2022109550} to genomics \cite{carstensen2010multivariate,gusto:hal-02682109}, ecology \cite{ward2022network}, epidemiology \cite{CHIANG2022505}, seismology  \cite{ogata1988statistical}, social sciences \cite{mohler2011self,farajtabar2016coevolve, crane2008robust} and stock prices moves \cite{Embrechts} 
to name but a few. 

Statistical inference for point processes has a long history \cite{Karr, Andersen}, with non-parametric approaches being favored for analyzing Poisson \cite{Kolaczyk, Patou, Willet-Nowak, MR3390131} and Hawkes processes \cite{Reynaud_Bouret_2010, Hansen:Reynaud:Rivoirard, MR4095329, Witten}. Maximum likelihood \cite{ogata1988statistical, MR961117, Anna} and the Bayesian framework have also shown good performance \cite{Rasmussen2013, HawkesBay1, HawkesBay2}. However, current frameworks are primarily based on the observation of a single point process, while increasingly, various fields are providing data in the form of replicated point processes. For example, in the field of earth sciences, point processes are essential to assess the seismic risk, which can have significant implications for civil engineering and insurance, for example. In the seminal work of \cite{ogata1988statistical}, dependent point processes were proposed to study preseismic quiescence and foreshocks that are expected to precede major earthquakes. This accurate modeling of earthquake hazards appears crucial to anticipate future earthquakes before they occur. Analyzing earthquake occurrences across various sites could help in identifying cities that show unusual earthquake accumulation patterns, or to more clearly identify the characteristic patterns of seismic activity before and after major earthquakes. In this work we will consider Turkish cities subject to recorded earthquakes over time, investigating the variability of these temporal patterns among cities \cite{Karasozen_etal_2018}. In single-cell genomics, it is now possible to individually characterize genomes based on epigenetic marks organized along chromosomes in a one-dimensional array \cite{zheng_massively_2017}. Furthermore, these repeated point patterns are instrumental in characterizing the diversity of cancer cells and their emergence of drug resistance. In neuroscience, observing spiking neurons has become routine for investigating whether their coordinated activity relates to the function of specific brain regions \cite{Cunningham2014}.  What unifies these examples is that replicated point processes offer an opportunity to study the variability of point patterns from a new perspective, characterizing data at a population level. This poses new mathematical, statistical, and methodological challenges, as most standard statistical frameworks are dedicated to continuous valued data and not to sets of points.

In this work we propose to put forward a different statistical perspective on point processes by developping a new framework for dimension reduction and visualisation of replicated point processes. Our framework is based on a functional analysis perspective -- to be contrasted to the perspective based on random sets. Viewing the realization of a point process as a random measure, and thus a bounded linear functional over an appropriately regular space of functions, we obtain uniform series representations and associated notions of principal modes of variation. This allows us to develop a Principal Component Analysis framework paralleling that of functional PCA (fPCA). fPCA is the workhorse of functional data analysis \cite{LWC13}. It enables the visualization of functional data, seen as random curves over compact sets (as described in \cite{silverman2002applied}), enables dimensionality reduction and regularisation in regression/testing/classification \cite{RO07,delaigle2012componentwise, hilgert2013minimax}, and provides a bridge to extend classical multivariate procedures (like model based clustering) to functional data \cite{EAV04,JP14}. The theoretical underpinnings of fPCA rely on two fundamental results: the Karhunen-Lo\`eve Theorem and  Mercer's Theorem \cite{hsing2015theoretical}. In particular, the Karhunen-Lo\`eve theorem asserts that the approximation of functional data by their principal components is uniformly convergent (in the mean-square sense). However, the theory of fPCA does not apply to point processes which are more intricate mathematical objects.

The link between point processes and functional data is a natural approach to define equivalent of PCA for samples of point processes, and has been attempted before, albeit via the intensity function. \cite{WMZ13} investigated PCA for doubly stochastic point processes in the temporal setting, while \cite{LG14} extended the analysis to spatio-temporal processes. \cite{BVAR06} also applied functional PCA to estimate the mean of a doubly stochastic Poisson process. These methods rely on kernel or projection estimators of the process intensities, leading to a Karhunen–Loève decomposition of the corresponding intensity functions. An alternative line of work, exemplified by \cite{BGMSS06} in the case of spatial processes, focuses on performing PCA on estimators of $L$-functions. However these approaches have notable drawbacks. First, instead of relying on the observed point process, they relate to an intensity (heuristically, a derivative) that is not observable. Furthermore, and most importantly, they do not allow for a Karhunen-Lo\`eve expansion for the measures associated to the initial process; this is a theoretical but also practical drawback, in terms of interpretation.  Other approaches include \cite{MYD07} who proposed a convenient Hilbertian setting for PCA of measures for the analysis of distribution functions appearing in the study of grain-size curve. Although the principle of the method is similar to ours, the setting is very different since the data is a sample of general random densities whose cumulative distribution function is observed on a fixed grid, rather than a sample of point processes.
 Closer in spirit to our approach, \cite{Carrizo22} recently obtained expansions of general random measures. More precisely, Theorem~4.1 of \cite{Carrizo22} proves that under mild assumptions, a finite regular random measure over $\R^d$ can be written as a series expansion of deterministic real finite measures weighted by uncorrelated real random variables with summable variances. The convergence of the series is in a weak sense, though, and hence does not constitute a Karhunen-Lo\`eve expansion. Consequently, the result is not used (or directly usable, for that matter) for the purposes of statistical analysis. 

To describe our contributions, we consider $n$ i.i.d. point processes $(N_1, \ldots, N_n)$,  and their associated random measures, $(\Pi_1,\ldots,\Pi_n)$. We start with Theorem~\ref{thm:KL}, which establishes the Karhunen-Lo\`eve expansion for the random measures $\Pi_i$, as well as a Mercer Theorem for the covariance measure of the process (see Theorem~\ref{thm:Mercer}).  We establish convergence in a strong sense. Utilizing the series expansion provided by the  Karhunen-Lo\`eve decompostion, we introduce the concept of principal measures, which can be seen as latent processes governing the dynamics of the observed measures. These results are based on the embedding of the measures $\Pi_i$ in a functional space using their associated \emph{cumulative mass functions} defined in Equation~\eqref{F-mu}. This embedding allows us to integrate our model within the functional PCA framework. We then construct a natural associated covariance operator whose eigenfunctions are the cumulative mass functions of the principal measures. Eigenelements are characterized for Poisson processes and Hawkes processes with exponential self-exciting function by solving second-order differential systems provided in Equation \eqref{sys-general}. When the Poisson process is homogeneous, the eigenfunctions of the covariance operator consist of a Fourier-like basis, and we also show that the eigenvalues exhibit a polynomial decay. When the Poisson process is inhomogeneous, the eigenfunctions are implicit; however, a detailed qualitative analysis allows us to show that the eigenfunctions are oscillatory. Surprisingly, we demonstrate that similar phenomena occur in the case of the Hawkes process. We propose a full estimation framework based on the observation of the $N_i$'s that does not require any smoothing, and provide convergence of estimators at the parametric rate. Our theoretical results are complemented by a simulation study that exemplifies the convergence of eigenvalues and eigenfunctions. We then apply our framework in seismology, genomics, and neuroscience, illustrating the versatility and power of our methodology. We also provide a R-package, \texttt{pppca}\footnote{\url{https://github.com/franckpicard/pppca}}, that implements our method.

The article is organized as follows. Section~\ref{sec:setting} introduces the setting of our work. Section~\ref{sec:KLMercer} establishes Karhunen-Lo\`eve and Mercer Theorems  for point processes. In Section~\ref{sec:illust}, we study the Poisson process and a specific class of Hawkes processes in detail by using the analytical point of view. Section~\ref{sec:RMPC} introduces our estimators of the eigenelements and establishes rates of convergence. Our numerical study is carried out in Section~\ref{sec:numerical_study}. Section~\ref{sec:applications} is devoted to applications. The proofs of our results are collected in the Supplementary Material.

\subsection*{Notation} 

We introduce some notation that will be useful throughout the article. The set $\mathbb N$ denotes the set of non-negative integers. For a set $I\subset\R^d$, for $d\geq 1$, we denote 
$$\mathbb L^2(I)=\Big\{ f:I\longmapsto\mathbb R:\ \|f\|<+\infty\Big\},$$
where $\|f\|$ stands for the $\mathbb L_2$-norm of $f$: $\|f\| = \left(\int_I |f(t)|^2dt\right)^{1/2}.$ The associated scalar product is denoted $\langle f,g\rangle = \int_I f(t)g(t)dt$. When $I=(0,1)$ and $\alpha\in\mathbb N$, we set $\partial_\alpha f = f^{(\alpha)}$. When $I=(0,1)^2$ and $\alpha=(\alpha_1,\alpha_2)\in\mathbb N^2$, we set $|\alpha|=\alpha_1+\alpha_2$ and $\partial_\alpha f = \frac{\partial^{|\alpha|} f}{\partial s^{\alpha_1}\partial t^{\alpha_2}}$.
 The derivatives and partial derivatives are taken in the distribution sense meaning that $\partial_{\alpha}f$ verifies the equality
	\[
	\langle \partial_\alpha f,\varphi\rangle  = (-1)^{|\alpha|}\langle f,\partial_{\alpha}\varphi\rangle, 
\] 
for all infinitely differentiable functions $\varphi$ with support in $I$.
\section{A signed-measure model for dimension reduction}\label{sec:setting}
Consider $n$ independent and identically distributed (i.i.d.) temporal point processes $(N_1, \ldots, N_n)$,  observed on the time interval $[0,1]$. These processes, defined on the probability space $(\Omega,\mathcal A,\mathbb P)$, are associated with some random measures $(\Pi_1,\hdots,\Pi_n)$ so that:
\[
\Pi_i(B)= \sum_{T \in N_i} \mathbf 1_{\{T\in B\}}  = {\rm Card}(N_i\cap B), \qquad B\in\mathcal B,
\]
where $\mathcal B$ is the Borel $\sigma$-field of subsets of $[0,1]$. Hence, $\Pi_1,\hdots,\Pi_n$ is a sample of random measures on $([0,1],\mathcal B)$. We denote by $\Pi$ a random measure following the same distribution as the $\Pi_i$'s and by $N$ the point process associated with $\Pi$. Assuming that $\mathbb E[\Pi([0,1])]<+\infty$, we denote the first moment of~$\Pi$ by $W$:
\[
W(B)=	\mathbb{E} \big[\Pi(B)\big],\qquad B\in\mathcal B.
\]
In the initial step of our dimension reduction framework, we center the observed data by defining a signed random measure 
$$\Delta_i(B) = \Pi_i(B)-W(B),\qquad B\in\mathcal B.$$
The $\Delta_i$'s have the same distribution as $\Delta=\Pi-W$.
Then, we define a second-order moment for $\Delta$, inspired from the formalism proposed by \cite{PZ16} and \cite{Carrizo22}. For this purpose, we define the signed measure of covariance: for all $B\times B'\in\mathcal B\otimes\mathcal B$
\begin{eqnarray*}
C_\Delta(B\times B')& :=& {\rm Cov}(\Pi(B),\Pi(B')) 
                               = \mathbb E \left[  \Delta(B)\Delta(B') \right] \\
                               &=& \mathbb E\Bigg[\sum_{T\in N,T'\in N}\mathbf 1_{\{(T,T')\in B\times B'\}}\Bigg]-( W \otimes W)(B\times B') 
\end{eqnarray*}
with $\mathcal B\otimes\mathcal B$ the product $\sigma$-algebra of $[0,1]\times[0,1]$, and where, for two measures $\mu$ and $\nu$ on $([0,1],\mathcal B)$ the product measure $\mu\otimes\nu$ on $([0,1]^2,\mathcal B\otimes \mathcal B)$ is defined such that \[
(\mu\otimes\nu)(B\times B') = \mu(B)\nu(B'), \qquad B,B'\in\mathcal B.
\]
We assume throughout that
\begin{equation}\label{eq:hyp2}
	\mathbb E\Big[\Pi^2([0,1])\Big]<\infty, 
\end{equation}
which implies that $C_\Delta$ is a finite measure. 
	
Now, to perform dimension reduction of the $N_i$'s (equivalently the  $\Pi_i$'s) we shall take advantage of the rich framework of the functional perspective. Following \cite{Cohn93} for instance, given a random measure $\mu$, we consider the associated \textit{cumulative mass function} defined by
\begin{equation}\label{F-mu}
F_\mu(t) = \mu( [0,t]), \quad t \in [0,1].
\end{equation}
Then, we introduce $K_{\Delta}$, the  bivariate cumulative mass function of the covariance measure $C_{\Delta}$, such that:
\[
K_{\Delta}(s,t)= C_{\Delta}([0,s]\times[0,t]) = \mathbb{E}[F_\Delta(s)F_\Delta(t)],\quad s,t\in[0,1]
\]  
and the associated integral operator
\[
\Gamma_{\Delta} = \mathbb{E} \Big( F_\Delta \otimes F_\Delta \Big)
\]
where $\otimes$ is the tensor product between two elements of $\mathbb L^2([0,1])$ denoted by $f\otimes g(h)=\langle f,h\rangle g$. The covariance operator $\Gamma_{\Delta}$ is expressed as
\begin{eqnarray*}
\Gamma_\Delta(f)(\cdot) &=&  \mathbb E (\langle f, F_\Delta \rangle F_\Delta(\cdot) ), \quad f \in {\mathbb L}^2, \\
&=& \int_0^1  K_{\Delta}(\cdot,t)f(t)dt.
\end{eqnarray*}
Then, dimension reduction is based on the Mercer expansion of the covariance kernel $K_{\Delta}$. Indeed, in the case where $K_\Delta$ is continuous, Mercer's Theorem applies and ensures the existence of a sequence of non-negative real numbers $(\lambda_j)_{j\geq 1}\in\ell^1(\mathbb N^*)$ and  orthonormal functional basis $(\eta_j)_{j\geq 1}$ of $\mathbb L^2([0,1])$ such that
\begin{equation}\label{Mercer}
K_\Delta(s,t)=\sum_{j\geq 1}\lambda_j\eta_j(s)\eta_j(t),\quad s,t\in[0,1],
\end{equation}
with uniform and absolute convergence of the series (see Theorem 4.6.5 of \cite{hsing2015theoretical}). Furthermore, $\eta_j$ is an eigenfunction of the operator $\Gamma_\Delta$ associated with the eigenvalue $\lambda_j$. 

Now, the question is the following. How can we use the decomposition~\eqref{Mercer} to build a series expansion for $\Pi$ and $C_\Delta$? This issue is tackled in the next section.

\section{Karhunen-Lo\`eve expansion and Mercer theorem for point processes}
\label{sec:KLMercer}

This section is devoted to prove a Karhunen-Lo\`eve theorem for the random measure $\Pi$ and a Mercer theorem for the covariance measure $C_\Delta$. The obtained expansions will be based on specific sequences of measures $(\mu_j)_j$, called \textit{principal measures}. 

We start from the decomposition~\eqref{Mercer}. Since $\Gamma_\Delta$ is a compact self-adjoint operator, the eigenvalues $(\lambda_j)_{j\geq 1} $ are isolated, with finite multiplicities and we have
$$\lambda_{j}\geq \lambda_{j+1}\geq 0,\quad\forall j\geq 1,$$
with $\lambda_j\to 0$ when $j\to+\infty$. We refer the reader to Sections~6.3 and 6.4 of \cite{Brezis11} for properties of the spectrum of a compact operator. In the sequel, we assume that in \eqref{Mercer}, 
\begin{equation}\label{lambdajpositive}
\lambda_j>0,\quad j=1,2,\ldots
\end{equation}
The following proposition is a central result that allows us to define our principal measures from the basis functions $(\eta_j)_{j\geq 1}$. 
\begin{proposition}\label{prop:muj} We suppose that Assumptions~\eqref{eq:hyp2} and \eqref{lambdajpositive} are satisfied and that $K_\Delta$ is continuous. Then, for all $j\geq 1$, the derivative in the distributional sense of $\eta_j$ is a measure, denoted by $\mu_j$, that verifies
\[
\eta_j(t)=\mu_j([0,t])=F_{\mu_j}(t),\quad t\in[0,1]
\]
and, for all  $\varphi\in\mathcal H_0^1=\Big\{f\in \mathbb L^2(I):f'\in \mathbb L^2(I) \text{ and }f(t) = 0 \text{ for all } t\notin I \Big\}$, 
\begin{equation}\label{phimu*}
\langle\eta_j,\varphi'\rangle =\langle F_{\mu_j},\varphi'\rangle =-\int_0^1\varphi(s)d\mu_j(s)=-\langle \varphi,\mu_j\rangle. 
\end{equation}
\end{proposition}
Now, using the $\mu_j$'s derived in Proposition~\ref{prop:muj} and following the proof of  Theorem 4.1 of \cite{Carrizo22},  we obtain that for all measurable and bounded function $\varphi$ on $[0,1]$, 
\[
\sum_{T\in N}\varphi(T)= \mathbb E\left[\sum_{T\in N}\varphi(T)\right]+\sum_{j\geq 1}\sqrt\lambda_j\xi_j\langle\mu_j,\varphi\rangle,
\]
where $\{\xi_j\}_{j\geq 1}$ is a sequence of uncorrelated real random variables of mean zero and variance one.
This result allows us to interpret the variations of $\sum_{T\in N}\varphi(T)$ around its mean as a sum of random variables with decreasing variance. It can also be expressed as follows: 
\begin{equation}\label{eq:KL}
\Pi(B)  = W(B) + \sum_{j\geq 1}\sqrt\lambda_j\xi_j\mu_j(B), \qquad B\in\mathcal B, 
\end{equation}
which implies for $t\in\mathcal [0,1]$,
\begin{equation}\label{eq:KL2}
F_{\Pi}(t) = \Pi([0,t])  = W([0,t]) + \sum_{j\geq 1}\sqrt\lambda_j\xi_j\mu_j([0,t]) = F_W(t)+ \sum_{j\geq 1}\sqrt\lambda_j\xi_jF_{\mu_j}(t). 
\end{equation}
The last equation means in particular that, under Assumption~\eqref{lambdajpositive}, the random sequence $(\xi_j)_{j\geq 1}$ is uniquely defined by the following equation 
\begin{equation}\label{def:scores1}
\xi_j = \frac{\langle F_\Pi-F_W,F_{\mu_j}\rangle}{\sqrt\lambda_j}, \qquad j\geq 1.
\end{equation}
To deduce the last expression, we have used that the functions $(F_{\mu_j})_{j\geq 1}$ are orthonormal. The convergence of the series holds then pointwise, i.e. for all $B\in\mathcal B$ (or for all measurable and bounded function $\varphi$). However, the strength of the Karhunen-Lo\`eve theorem for general stochastic processes is that the convergence holds uniformly on $[0,1]$ which prevents problems on the boundaries (typically Gibbs effect type phenomena). This is not ensured by Theorem 4.1 of  \cite{Carrizo22}. Therefore, we first establish a Karhunen-Lo\`eve decomposition for point processes for which the convergence holds uniformly for all sufficiently regular functions. Uniformity, however, cannot be in the total variation distance: as illustrated in Section~\ref{sec:illust} for Poisson and Hawkes processes, the measures $(\mu_j)_j$ are signed measures with very regular densities (typically $\mathcal C^\infty$). On the contrary, the measure $\Pi$ is discrete. This implies that the convergence cannot hold uniformly for all continuous and bounded functions, that is to say, uniformly in total variation. For this reason, in the sequel, uniform convergence will be related to Sobolev spaces instead. For this purpose, we introduce, for $I=(0,1)$ or $I=(0,1)^2$
	\[
	\mathcal H_0^{k} = \Big\{f\in \mathbb L^2(I):\partial_\alpha f\in \mathbb L^2(I) \text{ for all }|\alpha|\leq k\text{ and }f(t) = 0 \text{ for all } t\notin I \Big\}.
	\]

	The derivatives and partial derivatives are taken in the distribution sense.	 
	Sobolev norms of negative order are also introduced as the norm of the dual space of $\mathcal H_0^{k}$, namely
	 \[
\|\mu\|_{\mathcal H^{-k}} = \sup\Bigg\{|\langle f,\mu\rangle|:f\in\mathcal H_0^{k} \text{ and}\sum_{|\alpha|\leq k}\|\partial_{\alpha}f\|^2\leq 1 \Bigg\}. 
\]
The definition and properties of these spaces are described in Brezis \cite[Chapters 8, 9]{Brezis11}. With these definitions in place, we may now state:
\begin{theorem}[Karhunen-Lo\`eve Theorem for point processes]
\label{thm:KL}
Suppose that Assumptions~\eqref{eq:hyp2} and \eqref{lambdajpositive} are satisfied and that $K_\Delta$ is continuous. Then, there exists a sequence $\{\xi_j\}_{j\geq 1}$ of uncorrelated real random variables of mean zero and variance one such that
\[
\lim_{J\to+\infty}\mathbb E\left[\left\|\Pi-W-\sum_{j=1}^J\sqrt\lambda_j\xi_j\mu_j\right\|_{\mathcal H^{-1}}^2\right]=0.
\]
\end{theorem}
\begin{remark}
Using the fact that, for any Borel function $\varphi$,
$$
\langle \varphi,\Pi-W\rangle=\sum_{T\in N}\varphi(T)-\mathbb E\left[\sum_{T\in N}\varphi(T)\right]
$$
a direct consequence of Theorem~\ref{thm:KL} is that 
\begin{equation*}
\lim_{J\to+\infty}\sup_{\varphi\in \mathcal H^1_0,\, \|\varphi\|^2+ \|\varphi'\|^2\leq 1}\mathbb E\left[\left(\sum_{T\in N}\varphi(T)-\mathbb E\left[\sum_{T\in N}\varphi(T)\right]-\sum_{j=1}^J\sqrt{\lambda_j}\xi_j \left\langle\mu_j,\varphi\right\rangle\right)^2\right]=0.
\end{equation*}
\end{remark}
By analogy with fPCA, the elements of the sequence $\{\xi_j\}_{j\geq 1}$ are called \textit{the scores associated to the measure $\Pi$ (or associated to  the point process $N$)}. They play a central role in the study of the data we consider (see Section~\ref{sec:applications}). 
Indeed, for each individual $i$, it can be deduced from the Karhunen-Lo\`eve decomposition that for $J$ large, $\Pi_i$ is close to $W+ \sum_{j=1}^J\sqrt{\lambda_j}\xi_{i,j}\mu_j$, which can be written informally as
\begin{equation}\label{representation}
\Pi_i(B) \approx W(B)+ \sum_{j=1}^J\sqrt{\lambda_j}\xi_{i,j}\mu_j(B), \qquad B\in\mathcal B
\end{equation}
or
\begin{equation}\Pi_i([0,t]) \approx W([0,t])+ \sum_{j=1}^J\sqrt{\lambda_j}\xi_{i,j}F_{\mu_j}(t), \qquad t\geq 0
\end{equation}
with
\begin{equation}\label{def:scores}
\xi_{i,j} = \frac{\langle F_{\Delta_i},F_{\mu_j}\rangle}{\sqrt\lambda_j},
\end{equation}
by using \eqref{def:scores1} and $F_{\Delta_i}=F_{\Pi_i}-F_W$. 
As to what concerns the covariance measure $C_\Delta$, we establish the following Mercer theorem: 
\begin{theorem}[Mercer's Theorem for $C_\Delta$]
\label{thm:Mercer}
Suppose that Assumptions~\eqref{eq:hyp2} and \eqref{lambdajpositive} are satisfied and that $K_\Delta$ is continuous. Then
\begin{equation}\label{eq:CVunif}
\left\|\sum_{j=1}^J\lambda_j\mu_j\otimes\mu_j-C_{\Delta}\right\|_{\mathcal H^{-2}}\xrightarrow[J\to+\infty]{}0.
\end{equation}
\end{theorem}
In the next sections, we focus on the study of the eigenelements $(\lambda_j,F_{\mu_j})_{j\geq 1}$. Recall that we have an easy correspondence between the measures $\mu_j$ and the mass functions $F_{\mu_j}$ thanks to Eq.~\eqref{F-mu}. In Section~\ref{sec:illust}, we study Poisson and specific classes of Hawkes processes from the analytical point of view. For these central examples of point processes,  we derive the asymptotic behavior of $(\lambda_j,F_{\mu_j})$ when $j$ goes to $+\infty$. We even obtain the explicit form of the eigenelements for the homogeneous Poisson process for all $j$'s. In Section~\ref{sec:RMPC}, for the general setting, we provide an estimation scheme for the eigenelements based on realizations of $(N_1, N_2, ..., N_n)$. 
\section{Study of eigenelements for Poisson and Hawkes processes}\label{sec:illust}
The goal of this section is to complement the general results of Section~\ref{sec:KLMercer} with some concrete and explicit cases, by adopting a constructive point of view to determine how eigenelements  $(\lambda_j,F_{\mu_j})$ behave, at least asymptotically when $j$ goes to $+\infty$. We consider two classical point process models, namely Poisson processes and Hawkes  processes with exponential self-exciting functions.

\indent To present our results, we introduce, for $j\geq 1$, 
\begin{equation}\label{Hj}
H_j(t)=\int_t^1F_{\mu_j}(s)ds, \quad t\in [0,1]
\end{equation}
and consider solutions $(\lambda,y)$  of the following system of second-order integro-differential equations
\begin{equation}\label{sys-general}
\left\{
\begin{array}{ll}
\displaystyle -\lambda y''(t)=w(t)y(t)+\int_0^1 M(s,t)y(s)ds ,   &t\in(0,1),\\
y(1)  =0,\quad y'(0)=0,&
\end{array}
\right.
\end{equation}
where $\lambda>0$ and $y$ is twice continuously differentiable and 
\begin{enumerate}
\item for the Poisson process, $w$ is the intensity of the process and $M=0$,
\item for the Hawkes process with exponential self-exciting function, $w$ is a constant proportional to the baseline intensity and $M$ is an exponential convolution kernel (see Eq.~\eqref{SL-Hawkes}).
\end{enumerate}
We shall prove that the eigenelements $(\lambda_j,F_{\mu_j})_{j\geq 1}$ of $\Gamma_\Delta$ correspond to solutions $(\lambda_j,H_j)_{j\geq 1}$ of System \eqref{sys-general}, once $H_j$ is renormalized by a constant so that $\|F_{\mu_j}\|_2=1$ (observe that for any constant $c>0$, $(\lambda,y)$ is a solution of \eqref{sys-general} if and only if $(\lambda,c\times y)$ is also a solution). This renormalization is required since $(F_{\mu_j})_{j\geq 1}=(\eta_j)_{j\geq 1}$ is an orthonormal functional basis  of $\mathbb L^2([0,1])$.

Relationships connecting Karhunen-Lo\`eve expansions to differential equations of order 2 have already been established in the study of stochastic processes. See for instance Section A of \cite{CP15} for the case of the Ornstein-Uhlenbeck process. The system of equations \eqref{sys-general} appears naturally when we study the equation
\[
\lambda_j F_{\mu_j}(t)=\int_0^1K_\Delta(t,s)F_{\mu_j}(s)ds,\quad t\in[0,1]
\]
by taking derivatives of left and right hand sides. See for instance the straighfroward proof of Theorem~\ref{KL-SLP} in the supplementary file. In the sequel, except in Remarks~\ref{eigenvalue=0Poisson} and \ref{eigenvalue=0Hawkes}, we assume that the $\lambda_j$'s are positive.

\subsection{Eigendecomposition for Poisson point processes}
\label{subsec:Poisson_theory}

Here we suppose that $N$ is a Poisson process with intensity $t\in[0,1]\longmapsto  w(t)$ assumed to be continuous and positive on $(0,1)$. In this case, 
\[
W(B)=\int_B w(t)dt,\quad B\in\mathcal B
\]
and $\Pi(B)$ is a Poisson variable with parameter $W(B)$. Furthermore, for two disjoint Borel sets $B$ and $B'$, $\Pi(B)$ and $\Pi(B')$ are independent so that
\[
F_W(t)=W([0,t])=\int_0^t w(u)du, \quad 0\leq t\leq 1
\]
and
\[K_\Delta(s,t)=F_W\big(\min\{s;t\}\big)=\int_0^{\min\{s;t\}}w(u)du, \quad 0\leq s, t\leq 1.
\]
The following result, which involves the functions $H_j$ introduced in \eqref{Hj}, provides the PCA dimension reduction for Poisson processes.
\begin{theorem}\label{KL-SLP}
Assuming $w$ is continuous and positive on $(0,1)$, then $(\lambda_j,F_{\mu_j})_{j\geq 1}$ are the eigenelements of the operator $\Gamma_\Delta$ if and only if $(\lambda_j,H_j)_{j\geq 1}$ are solutions of the following Sturm-Liouville problem:
\begin{equation}\label{SLP}
\left\{
\begin{array}{ll}
 -\lambda y''(t)=w(t)y(t),   &t\in(0,1),\\
y(1)  =0,\quad y'(0)=0.&
\end{array}
\right.
\end{equation}
\end{theorem}
Mercer's and Karhunen-Lo\`eve expansions for Poisson processes are then provided by solutions of System \eqref{SLP} and by using the relationship
\[
F_{\mu_j}(t)=-H_j'(t),\quad t\in [0,1].
\]
\begin{remark}\label{eigenvalue=0Poisson}
The proof of Theorem~\ref{KL-SLP} shows that if $F_{\mu_0}$ is a continuous eigenfunction of $\Gamma$ associated with the zero eigenvalue, then $F_{\mu_0}$ is solution of \eqref{SLP} with $\lambda=0$ and then $F_{\mu_0}=0$. 
\end{remark}
Unfortunately, solutions of System \eqref{SLP} are typically not explicit. However, the specific case where $w$ is constant, corresponding to the homogeneous Poisson process, is of interest and can be easily dealt with.

\subsubsection{Homogeneous Poisson processes}\label{sec:HomogeneousPoisson}

In the case of homogeneous Poisson process, the function $w$ is constant on $[0,1]$: for some $w_0\in\R_+^*$,
\[
w(t)=w_0,\quad t\in [0,1]
\]
and the solutions of Equation~\eqref{SLP} can be obtained explicitly 

\[
\lambda_j=\frac{4w_0}{\pi^2(2j-1)^2},\quad H_j(t)=A_j\cos(\pi(2j-1)t/2), \quad t\in[0,1],\quad j\geq 1,\]
where, in the last expression, $A_j$ can be any constant independent of $t$. But, taking into account the constraint $\|H_j'\|=\|F_{\mu_j}\|=1$, we must have $|A_j|=\sqrt{2}/(\pi j)$. Without loss of generality, we take $A_j=-\sqrt{2}/(\pi j)$ so that
\begin{equation}\label{sin-Fj}
F_{\mu_j}(t)= \sqrt{2}\sin(\pi(2j-1)t/2),\quad t\in [0,1] ,\quad j\geq 1. 
\end{equation}

In particular, the eigenfunctions $F_{\mu_j}$ are highly oscillating functions with exactly $j$ zeros on $[0,1]$ at points $x_k=2k/(2j-1),$ $k=0,\ldots,j-1$.

\bigskip
Then, we can finally write that if $N$ is a homogeneous Poisson process with intensity $w_0$,
\[
\Pi([0,t]) = w_0t +\frac{2\sqrt{2w_0}}{\pi}\sum_{j\geq 1}\frac{\xi_j}{2j-1}\sin(\pi(2j-1)t/2),\quad t\in [0,1],
\]
with $\{\xi_j\}_{j\geq 1}$ uncorrelated centered real random variables of unit variance.
\begin{remark}
Observe that in the homogenous Poisson setting, $K_\Delta(s,t)=w_0\min\{s;t\}$, which corresponds to the covariance kernel of the Brownian motion, so that 
 the Karhunen-Lo\`eve decomposition of $\Pi$ can be deduced from the one of the Brownian motion (see e.g. \cite{ash_topics_1975}, pp. 41--42). The main difference with the Brownian motion is that the scores $\{\xi_j\}_{j\geq 1}$ are not Gaussian. Determining the distribution of the scores in the Poisson setting remains an open and interesting problem. It is beyond the scope of this paper.
\end{remark}

\subsubsection{Inhomogeneous Poisson processes}
\label{subsubsec:inhPP}

The general Poisson case is much more involved, since, when $w$ is not constant, we cannot derive explicit solutions of \eqref{SLP}, in full generality. But we prove in this section that solutions of \eqref{SLP} exist and some useful qualitative properties of them can be obtained. Actually, system \eqref{SLP} is a self-adjoint Sturm-Liouville problem with separated boundary conditions, which has been extensively considered in the literature. We refer the reader to Chapter 4 of \cite{Zettl-book} and more precisely to Equations (4.1.1)-(4.1.4) with
$$p\equiv1,\quad q\equiv0,\quad \lambda^{-1} \mbox{ instead of } \lambda.$$
Therefore, we have the following result \cite{Zettl-book}.
\begin{theorem}[Theorems 4.3.1 and 4.6.2 of \cite{Zettl-book}]\label{thm:inhPP}
We consider the solutions $(\lambda_j,H_j)_{j}$  of the Sturm-Liouville problem \eqref{SLP} with $w\in{\mathbb L}_1[0,1]$ and assume that $w$ is positive on $(0,1)$. We have:
\begin{enumerate}
\item The eigenvalues $(\lambda_j)_{j\geq 1} $ are simple and asymptotically, we have
$$\lambda_j  \underset{j\to+\infty}{\sim}\frac{\Big(\int_0^1\sqrt{w(u)}du\Big)^2}{(j\pi-\pi/2)^2.}$$
\item For $j\geq 1$, if $H_j$ is the solution associated with $\lambda_j$, then $H_j$ is unique up to a multiplicative constant and has exactly $j-1$ zeros on $(0,1)$.
\item The sequence $(H_j)_{j\geq 1}$ can be normalized to be an orthonormal sequence in ${\mathbb L}_2([0,1],w)$, i.e.
$$
\int_0^1H_j(u)H_{j'}(u)w(u)du=\left\{
\begin{array}{ccc}
0  &\mbox{if}   &j\not=j',   \\
 1 & \mbox{if}   &j=j'.   \\
 \end{array}
\right.
$$
\end{enumerate}
\end{theorem}

\begin{remark}
By using the third point of Theorem~\ref{thm:inhPP}, we observe that for $j\not=j'$,
\begin{align*}
\int_0^1F_{\mu_j}(u)F_{\mu_{j'}}(u)du&=\int_0^1H'_j(u)H'_{j'}(u)du\\&=H'_j(1)H_{j'}(1)-H'_j(0)H_{j'}(0)-\int_0^1H''_j(u)H_{j'}(u)du\\
&=\lambda_j^{-1}\int_0^1H_j(u)H_{j'}(u)w(u)du=0,
\end{align*}
so that, once conveniently normalized, $(F_{\mu_j})_{j\geq 1}$ is an ${\mathbb L}_2$-orthonormal basis, as required by our theory. 
\end{remark}
Unlike the homogeneous case, we do not obtain the exact values of the eigenvalues $\lambda_j$ for all $j$'s, but the first point of Theorem~\ref{thm:inhPP} provides their asymptotic behavior when $j\to+\infty$. The numerical study of Section~\ref{sec:numerical_study} goes further since it shows that, for even small values of $j$, $\lambda_j$ is close to $(\int_0^1 \sqrt{w(u)}du)^2/(j \pi - \pi/2)^{2}$, exactly as for the homogeneous case (see Fig.~\ref{Fig:PP_eigenval}).

Combining the second point of Theorem~\ref{thm:inhPP} and the boundary condition $H_j(1)=0$, we obtain that $H_j$ has exactly $j$ zeros on $(0,1]$. Since $F_{\mu_j}(t)=-H_j'(t)$ for $t\in [0,1]$ and $F_{\mu_j}(0)=0$, Rolle's theorem yields that $F_{\mu_j}$ has at least $j$ zeros on $[0,1)$ and is therefore a highly oscillating functions, exactly as for the homogeneous case for which the functions $F_{\mu_j}$ are sine functions (see \eqref{sin-Fj}). These qualitative properties are confirmed by the numerical study of Section~\ref{sec:numerical_study} (see Fig.  \ref{Fig:PP_eigenfun}). The numerical examples considered show that the $F_{\mu_j}$'s have exactly $j$ zeros on $[0,1)$.

To conclude, when $w$ is not constant, we cannot give a closed expression of eigenelements of the operator $\Gamma_\Delta$, but through the study of their behavior provided by the Surm-Liouville theory, we have established that they can be viewed as perturbations of the eigenelements of the homogeneous Poisson case.

\subsection{Eigendecomposition for Hawkes processes with exponential self-exciting functions}\label{sec:Hawkes}

We extend our analytical results by providing results on dimension reduction for (linear) Hawkes processes. Hawkes processes were designed by Alan G. Hawkes \cite{Hawkes71a, Hawkes71b} to model the properties of self-excitation. This model seems particularly natural to account for dependencies in our point process model. Hawkes processes can be viewed as a branching process over a homogeneous Poisson process but they are also characterized by their stochastic intensity $w(t)$ defined by 
\[w(t)=\lim_{\delta\to 0}\delta^{-1}\mathbb E\Big[\Pi[0,t+\delta]-\Pi[0,t]|{\mathcal F}_t\Big],\]
where the filtration ${\mathcal F}_t$ stands for the information available up to (but not including) time $t$; $w(t)$ can also be interpreted as the probability to have a new occurence at time $t$ given the past.
\begin{definition}	
The (linear) Hawkes process with baseline intensity $w_0>0$ and self-exciting function $h:\R_+\longmapsto\R_+$ is a point process $N$ with intensity
\begin{equation}\label{H-intensity}
w(t)=w_0+\int_{-\infty}^{t^-}h(t-s)dN_s=w_0+\sum_{T\in N, \ T<t}h(t-T).
\end{equation}
\end{definition}
Observe that the larger $h$, the stronger the influence of the past. If $h=0$, the associated Hawkes process corresponds to a homogeneous Poisson process whose intensity is the baseline intensity $w_0$. In the sequel, we mainly consider an exponential self-exciting function, i.e.
\begin{equation}\label{h-exp}
h(t)=\alpha\exp(-\beta t),\quad t>0,
\end{equation}
and in this case, the intensity becomes
\[
w(t)=w_0+\alpha\int_{-\infty}^{t^-}\exp(-\beta (t-s))dN_s=w_0+\alpha\sum_{T\in N, \ T<t}\exp(-\beta (t-T)),
\]
meaning that the influence of the past decreases exponentially in time. The intensity associated with an exponential self-exciting function becomes a Markov process \cite{Oakes75}.

In the sequel, we assume that $\|h\|_1<1$. For the case of exponential self-exciting function $h$ given in \eqref{h-exp}, this condition is equivalent to $\alpha<\beta$.  
Assuming  $\|h\|_1<1$ ensures that  there exists a unique stationary version of the Hawkes process  $N$ with intensity given by \eqref{H-intensity} (see Theorem~1 of \cite{BM1}). Furthermore, under this assumption, the mean and the covariance kernel of $N$ are well-defined and if we consider the exponential case, we obtain explicit expressions:
\begin{proposition}\label{cor:Hawkes_exp}
We assume that for $t> 0$, $h(t)=\alpha e^{-\beta t}$, with $0<\alpha<\beta$. Then,
\begin{equation}\label{mean-Hawkes}
 F_W(t)=W([0,t])=\frac{w_0\beta t}{\beta-\alpha}, \quad t>0,
\end{equation}
and for $0<s\leq t$,
\begin{equation}\label{eq:K_Hawkes_exp}
K_\Delta(s,t)=\frac{\beta^3w_0}{(\beta-\alpha)^3}s+\frac{w_0\alpha\beta(2\beta-\alpha)}{2(\beta-\alpha)^4}\big(-1-e^{(\alpha-\beta)(t-s)}+e^{(\alpha-\beta)t}+e^{(\alpha-\beta)s}\big).
\end{equation}
\end{proposition}
The proof of Proposition~\ref{cor:Hawkes_exp} can easily be deduced from results in Section 3 of \cite{GaoZhu2018}.
Subsequently, using the functions $H_j$ introduced in \eqref{Hj}, the following result provides the PCA decomposition for Hawkes processes with exponential self-exciting function when studied on the compact interval $[0,1]$.
\begin{theorem}\label{Hawkes-ed}
We assume that for $t> 0$, $h(t)=\alpha e^{-\beta t}$, with $0<\alpha<\beta$. Then, $(\lambda_j,F_{\mu_j})_{j\geq 1}$ are the eigenelements of the operator $\Gamma_\Delta$ if and only if $(\lambda_j,H_j)_{j\geq 1}$ are solutions of the following system:
\begin{equation}\label{SL-Hawkes}
\left\{
\begin{array}{ll}
\displaystyle -\lambda y''(t)=\frac{\beta w_0}{\beta-\alpha}y(t)+\frac{w_0\alpha\beta(2\beta-\alpha)}{2(\beta-\alpha)^2}\int_0^1e^{-(\beta-\alpha)|t-s|}y(s)ds,   &\quad t\in(0,1),\\
y(1)  =0,\quad y'(0)=0.&
\end{array}
\right.
\end{equation}
\end{theorem}
\begin{remark}\label{eigenvalue=0Hawkes}
Similarly to the Poisson case, if $F_{\mu_0}$ is a continuous eigenfunction of $\Gamma$ associated with the zero eigenvalue, then $F_{\mu_0}=0$.
\end{remark}
The following result shows that under Assumption~\eqref{cond-rho}, the solutions of the ordinary differential equation \eqref{SL-Hawkes} are oscillating functions. In this particular setting, we then obtain a behavior close to the behavior observed for the Poisson case. 
\begin{theorem}\label{Hawkes-ed-sol}
Assume $\alpha<\beta$ and
\begin{equation}\label{cond-rho}
\frac{\alpha(2\beta-\alpha)}{\beta-\alpha}\Big(2+\frac{3-2e^{-(\beta-\alpha)/2}-e^{-(\beta-\alpha)}}{\beta-\alpha}\Big)<1.
\end{equation}
Let \[w_1=\frac{\beta w_0}{\beta-\alpha}.\] 
Then, the eigenvalues $(\lambda_j)_{j\geq 1} $ are simple and for all $j\geq 1$, there exists a solution $(\lambda_j,H_j)$ to the system~\eqref{SL-Hawkes} such that $\|F_{\mu_j}\|_2=1$ and we have
\begin{equation}\label{lambdaj}
\Bigg|\lambda_j-\frac{4w_1}{\pi^2\big(2j-1\big)^2}\Bigg|\leq C_1j^{-4}
\end{equation}
and
\begin{equation}\label{etaj}
\sup_{t\in[0,1]}\big|F_{\mu_j}(t)- \sqrt{2}\sin(\pi(2j-1)t/2)\big|\leq C_2j^{-1},
\end{equation}
for $C_1$ and $C_2$ two constants not depending on $j$. 
\end{theorem}
Let us discuss Assumption~\eqref{cond-rho} which is used to to solve the system~\eqref{SL-Hawkes}. More precisely, this condition allows to control norms of operators introduced in Lemma~\ref{Opsup} (in the Supplementary File) to obtain the expression of $H_j$ provided by Equation~\eqref{Fjexpan}. Whether Assumption~\eqref{cond-rho} is necessary to state the conclusions of Theorem~\ref{Hawkes-ed-sol} remains an open question. Since $$2+\frac{3-2e^{-(\beta-\alpha)/2}-e^{-(\beta-\alpha)}}{\beta-\alpha}\in [2,4]$$ for any $0<\alpha<\beta$, the main term in the left hand side of Assumption~\eqref{cond-rho} is $$\frac{\alpha(2\beta-\alpha)}{\beta-\alpha}=\frac{\alpha(2-\|h\|_1)}{1-\|h\|_1}.$$ It means that Assumption~\eqref{cond-rho} is satisfied as soon as $\alpha$ is not too large and $\|h\|_1$ is not too close to 1. Interpreting $\alpha$ as the size of the jumps and $\|h\|_1$ as the influence of the past, it means that the size of the jumps cannot be too large and the influence of the past cannot be too strong. Note that when $\|h\|_1$ or  $\alpha$ go to 0, the Hawkes process becomes closer to a homogeneous Poisson process with parameter $w_1$. This can be seen in the definition of the intensity $w$ but also 
in the reformulation of System~\eqref{SL-Hawkes}:
$$
\left\{
\begin{array}{ll}
\displaystyle -\lambda y''(t)=w_1 y(t)+\frac{w_1\alpha(2-\|h\|_1)}{2(1-\|h\|_1)}\int_0^1e^{-\alpha\big(\frac{1}{\|h\|_1}-1\big)|t-s|}y(s)ds,   &\quad t\in(0,1),\\
y(1)  =0,\quad y'(0)=0.&
\end{array}
\right.
$$
If $y$ is bounded, then the second term of the right hand side goes to 0 when $\|h\|_1$ or  $\alpha$ goes to~0.

Surprisingly, we observe that the eigenvalues have an asymptotic behavior similar to the homogeneous Poisson case. Since the upper bound \eqref{lambdaj} is polynomial in $j$, approximations are sharp. Observe that the parameter $w_1$ is the average number of points in $[0, 1]$. It is the analog of parameter $w_0$ for the homogeneous Poisson process studied in Section~\ref{sec:HomogeneousPoisson}. In both cases, we have that $\lambda_j$ is close to $\frac{{\mathbb E}[\Pi([0,1])]}{(j\pi-\pi/2)^2}.$ The Poisson analogy for eigenfunctions $F_{\mu_j}$ is valid as well (see \eqref{sin-Fj}). To complement these theoretical contributions, a numerical study is lead in Section~\ref{sec:numerical_study}. In particular, we examine in Section~\ref{sec:numerical_study} whether Assumption~\eqref{cond-rho} is a technical artefact or not.

\section{Estimation of eigenelements}\label{sec:RMPC}

In practice, dimension reduction of point processes is performed by estimating the principal measures $\mu_j$, the associated eigenvalues $\lambda_j$ and the scores   $\xi_{i,j}$ (see \eqref{representation}). The estimation of the principal measures $\mu_j$ is based on the estimation of the associated cumulative mass functions $F_{\mu_j}$ and hence necessitates the estimation of eigenfunctions and eigenvalues of the covariance operator $\Gamma_{\Delta}$. To proceed, by using the random measures $(\Pi_1,\ldots,\Pi_n)$, we first define an empirical version of the covariance operator $\Gamma_\Delta$ using the empirical counterparts of the covariance measure $C_\Delta$ and kernel $K_\Delta$. We recall that we assume that the point processes $(N_1, \ldots, N_n)$ are independent and identically distributed.
Let  
$$
\widehat{W} = \frac{1}{n} \sum_{i=1}^n \Pi_i \quad\text{ and }\quad\widehat{\Delta}_i = \Pi_i-\widehat{W}.
$$
Then we define 
\[
\widehat C_{\widehat\Delta}(B \times B')=\frac1n\sum_{i=1}^n\sum_{T, T'\in N_i}\mathbf 1_{\{(T,T')\in B \times B'\}}-\widehat W (B) \times\widehat W(B'), \quad B, B' \in\mathcal B,
\]
which provides the definition of the empirical covariance kernel:
\[
\widehat K_{\widehat\Delta}(s,t)=\widehat C_{\widehat\Delta}([0,s]\times[0,t]),\quad s,t\in[0,1]
\]  
along with the empirical integral operator:
\[
\widehat\Gamma_{\widehat\Delta}(f)(\cdot) = \int_0^1\widehat K_{\widehat\Delta}(\cdot,t)f(t)dt,
\]
that can be expressed using the cumulative mass functions of the empirical measures, namely $F_{\widehat\Delta_i}(t)=\widehat\Delta_i([0,t])$, which leads to 
\begin{equation*}\label{dec-hatGamma}
\widehat\Gamma_{\widehat\Delta}=\frac1n\sum_{i=1}^nF_{\widehat\Delta_i}\otimes F_{\widehat\Delta_i}.
\end{equation*}
 Then, ${\rm Im}(\widehat\Gamma_{\widehat\Delta})\in{\rm span}\{F_{\widehat\Delta_1},\hdots,F_{\widehat\Delta_n}\}$, which implies that $\widehat\Gamma_{\widehat\Delta}$ is a finite-rank operator. Since it is also self-adjoint, the covariance operator is compact and the diagonalization theorem ensures the existence of a basis $(\widehat\eta_j)_{j\geq 1}$ of eigenfunctions of $\widehat\Gamma_{\widehat\Delta}$. We denote by $(\widehat\lambda_j)_{j\geq 1}$ the associated eigenvalue sequence, sorted in non-increasing order, so that:
 $$
 \widehat\Gamma_{\widehat\Delta}= \sum_{j \geq 1} \widehat{\lambda}_j \widehat{\eta}_j \otimes \widehat{\eta}_j.
 $$
 In particular, the eigenelements $(\lambda_j,\eta_j)_{j\geq 1}$ are estimated by $(\widehat{\lambda}_j, \widehat{\eta}_j)_{j\geq 1}$.  Finally, for any process $N_i$ with $i\in\{1,\ldots,n\}$ and any $j\geq 1$, we estimate the scores  $\xi_{i,j}$ defined in Section~\ref{sec:KLMercer} by
\[
\widehat\xi_{i,j} = \frac{\langle \widehat \eta_j,F_{\widehat\Delta_i}\rangle}{\sqrt{\widehat\lambda_j}} \qquad\text{ if }\widehat\lambda_j>0. 
 \]
The case $\widehat\lambda_j=0$ is not considered since it corresponds to eigenfunctions $\widehat\eta_j$ that are in the kernel of $\widehat\Gamma_{\widehat\Delta}$; this is not encountered in practice. 
 
 Since the operator $\widehat\Gamma_{\widehat\Delta}$ is a finite-rank operator, it is possible to build a matrix, denoted by $\widehat G_{\widehat\Delta}$, such that the eigenelements of $\widehat\Gamma_{\widehat\Delta}$ can be derived explicitly from the eigenelements of $\widehat G_{\widehat\Delta}$.  To construct $\widehat G_{\widehat\Delta}$, we consider all occurrences sorted in non-decreasing order:
\[
\mathcal{T} = \Big(\bigcup_{i=1}^n N_i\Big)\cup\{0,1\}=\Big\{T_0,\hdots,T_{|\mathcal{T} |}\Big\}, 
\]
so that $T_0 = 0$ and $T_{|\mathcal{T} |}=1.$ If some ties occur when we take the union of the processes $N_i$, we identify them in only one $T_\ell$ so that for any $\ell\not=\ell'$, we have $T_\ell\not= T_{\ell'}$. Then, we define the histogram system $(e_1, \hdots, e_{|\mathcal{T} |})$  associated to this grid
\[
e_\ell(t)=\frac1{\sqrt{T_{\ell}-T_{\ell-1}}}1_{[T_{\ell-1};T_{\ell})}(t),\quad \ell=1,\hdots,|\mathcal{T} |,
\]
 and we set
 \[
\widehat G_{\widehat\Delta}=\left(\langle \widehat\Gamma_{\widehat\Delta} (e_\ell),e_{\ell'}\rangle\right)_{1\leq \ell,\ell'\leq |\mathcal{T} |}.
\]
Since $ \widehat\Gamma_{\widehat\Delta}$ is self-adjoint operator, $\widehat G_{\widehat\Delta}$ is a self-adjoint matrix. In the following Lemma, we show that $\widehat G_{\widehat\Delta}$ can be easily computed. Furthermore, we show that eigenelements of $ \widehat\Gamma_{\widehat\Delta}$ can be deduced from eigenelements of $ \widehat G_{\widehat\Delta}$.
\begin{lemma}\label{eq:construct_GDelta}
\text{}
\begin{enumerate}
\item The elements of $\widehat G_{\widehat\Delta}$ are in a closed form: for any $1\leq\ell,\ell'\leq |\mathcal{T} |$,
$$\left( \widehat G_{\widehat\Delta} \right)_{\ell, \ell'}=
			\frac{\sqrt{(T_\ell-T_{\ell-1})(T_{\ell'}-T_{\ell'-1})}}{n}  \times  \sum_{i=1}^n\Big(F_{\Pi_i}(T_{\ell'-1})F_{\Pi_i}(T_{\ell-1}) -F_{\widehat W} (T_{\ell'-1})F_{\widehat W}(T_{\ell-1})\Big).
$$
In particular, if for all $\ell$, there exists a unique $i$ such that $T_\ell\in N_i$, we have $F_{\widehat W}(T_{\ell-1})=(\ell-1)/n$ and
$$\left( \widehat G_{\widehat\Delta} \right)_{\ell, \ell'}=
			\frac{\sqrt{(T_\ell-T_{\ell-1})(T_{\ell'}-T_{\ell'-1})}}{n}  \times  \sum_{i=1}^n\left(F_{\Pi_i}(T_{\ell'-1})F_{\Pi_i}(T_{\ell-1}) -\frac{\ell-1}n\frac{\ell'-1}n\right).
$$
\item Let ${\rm Sp}(\widehat\Gamma_{\widehat\Delta})$ (resp. ${\rm Sp}(\widehat G_{\widehat\Delta})$) be the set of eigenvalues of the operator $\widehat\Gamma_{\widehat \Delta}$ (resp. of $\widehat G_{\widehat\Delta}$). We have:
\[
{\rm Sp}(\widehat\Gamma_{\widehat\Delta})\backslash\{0\}\subset {\rm Sp}(\widehat G_{\widehat\Delta})\subset {\rm Sp}(\widehat\Gamma_{\widehat\Delta}). 
\]
\item Let $\widehat v_j=(\widehat v_j^1,\hdots,\widehat v_j^{|\mathcal{T}|})^t$ a unit-norm eigenvector of $\widehat G_{\widehat\Delta}$ associated with the eigenvalue $\widehat\lambda_j$ and let
	$$
	\widehat\eta_j = \sum_{\ell=1}^{|\mathcal{T}|}\widehat v_j^\ell e_\ell. 
	$$
	Then $\widehat\eta_j$ is a unit-norm eigenfunction of the operator $\widehat\Gamma_{\widehat\Delta}$ associated with the eigenvalue $\widehat\lambda_j$ and there exists a unique measure $\widehat\mu_j$ such that  $F_{\widehat\mu_j} = \widehat\eta_j. $
\end{enumerate}
\end{lemma}
\begin{remark}
If $0\in{\rm Sp}(\widehat G_{\widehat\Delta})$, then the previous result shows that ${\rm Sp}(\widehat G_{\widehat\Delta})= {\rm Sp}(\widehat\Gamma_{\widehat\Delta}).$
\end{remark}
\begin{remark}
Of course, the definition of $\widehat\eta_j$ is arbitrary since several eigenvectors may be associated to the eigenvalue $\widehat\lambda_j$. In particular, we could consider $-\widehat\eta_j$. Actually, in practice, the dimension of the estimated eigenspaces associated to each $\widehat\lambda_j$ is equal to 1. In this case, the definition of $\widehat\eta_j$ is unique up to its sign.
\end{remark}
\begin{remark}
The estimated scores can be written
\begin{eqnarray*}
\widehat\xi_{i,j}&=&\frac{\langle \widehat\eta_j,F_{\widehat\Delta_i}\rangle}{\sqrt{\widehat\lambda_j}}=\widehat\lambda_j^{-1/2}\sum_{\ell=1}^{|\mathcal T|}\frac{\widehat v_j^\ell}{\sqrt{T_{\ell-1}-T_\ell}} \int_{T_{\ell-1}}^{T_\ell}F_{\widehat\Delta_i}(t)dt\\
&=&	\widehat\lambda_j^{-1/2}\sum_{\ell=1}^{|\mathcal{T}|}\widehat v_j^\ell \sqrt{T_\ell-T_{\ell-1}}F_{\widehat\Delta_i}(T_{\ell-1}),
\end{eqnarray*}

where the last equality comes from Equation~(S3.1.24) of the supplementary file.
\end{remark}
\begin{remark}
Since $\widehat\Gamma_{\widehat\Delta}$ is a finite-rank operator, we establish in Section~\ref{sec:proof:op} of the supplementary file that ${\rm Im}(\widehat\Gamma_{\widehat\Delta})\subset{\rm span}\{e_1,\hdots,e_{|\mathcal T|} \}$ (see \eqref{eq:incluImST}). In particular, eigenvectors associated with non-zero eigenvalues of $\widehat\Gamma_{\widehat\Delta}$ are histograms. The expression provided by the third point of Lemma~\ref{eq:construct_GDelta} specifies the values taken by each histogram $
	\widehat\eta_j$.
\end{remark}

The following theorem shows that the estimates introduced previously achieve the parametric rate for estimating eigenelements $\lambda_j$, $\eta_j$ and $\mu_j$. \begin{theorem}\label{thm:CV_rates}
We assume that 
\begin{eqnarray}\label{hyp:momdelta4}
\mathbb E[\|F_\Delta\|^4]<+\infty.
\end{eqnarray}
We have:
\begin{equation}\label{lambda-est}
\mathbb E\Big[\sup_{j\geq 1}|\widehat\lambda_j-\lambda_j|^2\Big]\leq 4\frac{\mathbb E[\|F_\Delta\|^4]}{n}.
\end{equation}	
If we further assume that the eigenvalues of $\,\Gamma_\Delta$ are simple, then, by introducing the eigengaps \[\delta_1=\lambda_1-\lambda_2\]
and
\[\delta_j=\min\left\{\lambda_j-\lambda_{j+1};\lambda_{j-1}-\lambda_{j}\right\},\quad j\geq 2,\] with $\widetilde\eta_j={\rm sign}(\langle\widehat\eta_j,\eta_j\rangle)\eta_j$, 
we have 
\begin{equation}\label{eta-est}
\mathbb E[\|\widehat\eta_j-\widetilde\eta_j\|^2]\leq32\delta_j^{-2}\frac{\mathbb E[\|F_\Delta\|^4]}{n}
\end{equation}
and
\begin{equation}\label{mu-est}
		\mathbb E\left[\|\widehat\mu_j-\widetilde\mu_j\|^2_{\mathcal H^{-1}}\right]\leq 32\delta_j^{-2}\frac{\mathbb E[\|F_\Delta\|^4]}{n},
		\end{equation}
with $\widetilde\mu_j = {\rm sign}(\langle\widehat\eta_j,\eta_j\rangle)\mu_j$. 		
\end{theorem}
The moment Assumption~\eqref{hyp:momdelta4} is very mild. Furthermore, Theorem~\ref{thm:inhPP} shows that for Poisson processes, the eigenvalues of $\Gamma_\Delta$ are simple. For Hawkes processes studied in Section~\ref{sec:Hawkes}, this condition holds under Assumption~\eqref{cond-rho} (see Theorem~\ref{Hawkes-ed-sol}).
Even if our framework is nonparametric in nature, let us observe that regularization is not necessary to achieve parametric rates.

\section{Numerical study}\label{sec:numerical_study}

We now complement our theoretical contributions with some empirical results. We probe, in particular, the following questions: 
 how does the truncated approximation behave? In particular, to what extent are the properties of eigenfunctions preserved in a non-asymptotic setting? For Hawkes processes, some restrictions due to Assumption~\eqref{cond-rho} appear. Are these limitations purely theoretical ones? We discuss these points in the sequel.

\begin{figure}
	\begin{center}
		\includegraphics[scale=0.5]{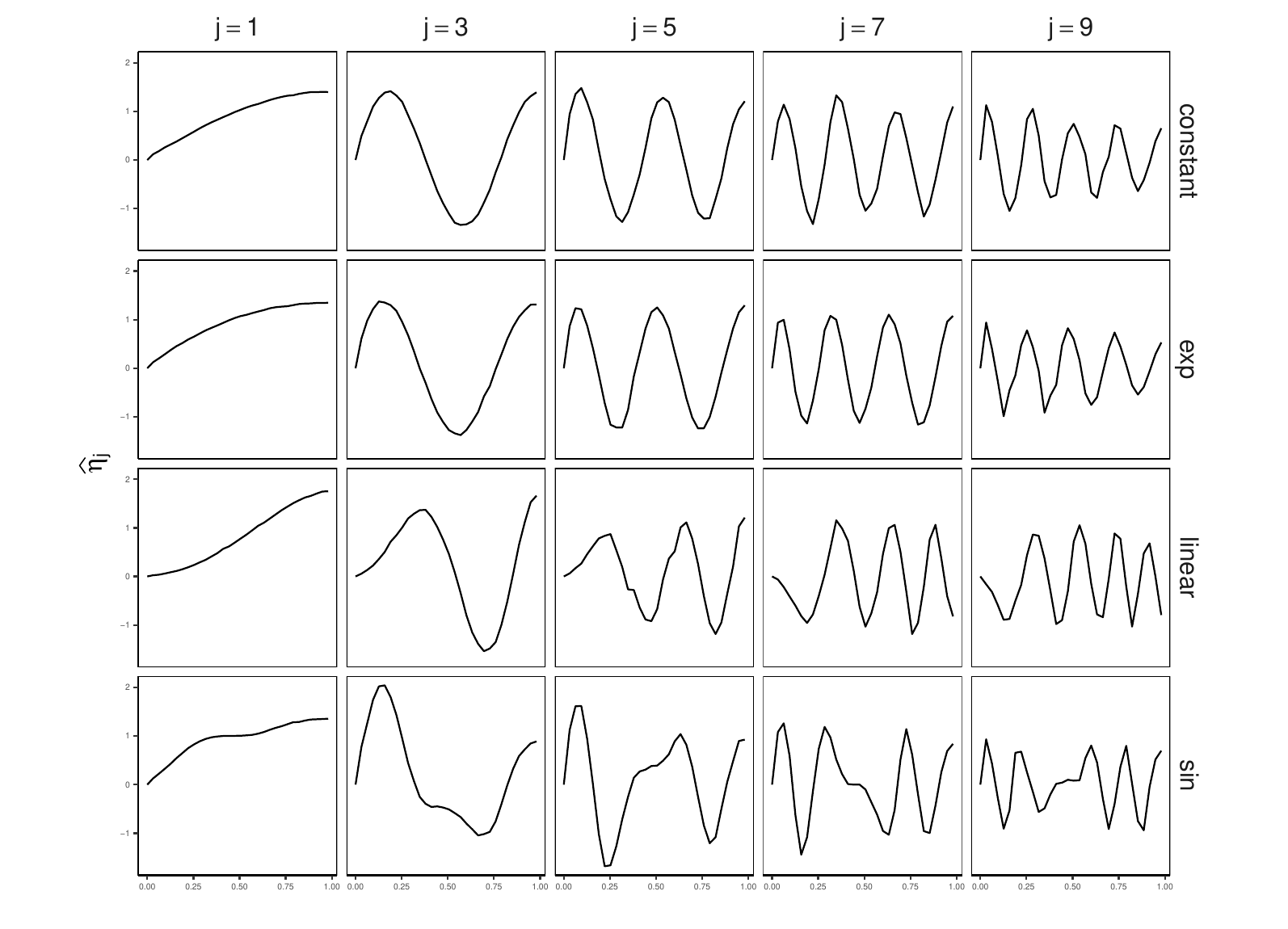}
	\end{center}
\caption{Average eigenfunctions for Poisson processes with different intensity functions over 50 replicates.} \label{Fig:PP_eigenfun}
\end{figure}

\begin{figure}
	\begin{center}
		\includegraphics[scale=0.8]{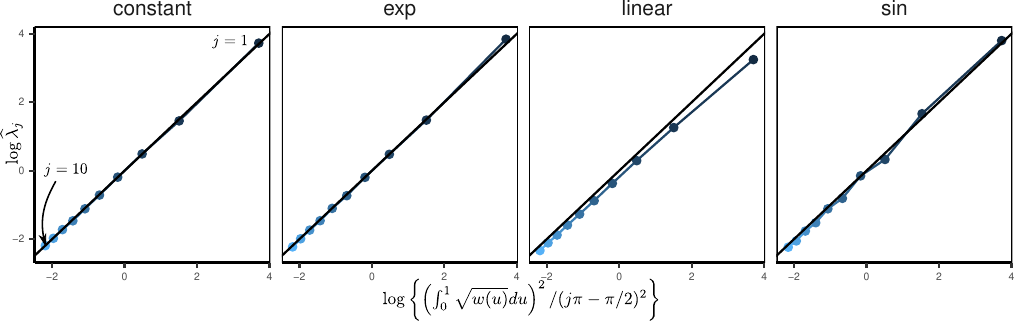}
	\end{center}
\caption{
	Average eigenvalues (log-scale) for Poisson Processes over 50 replicated. Each dot corresponds to a value of $j \in \{1,\hdots, 10\}$. The empirical average is plotted vs the expected theoretical asymptotic regime of eigenvalues in $(\int_0^1 \sqrt{w(u)}du)^2/(j \pi - \pi/2)^{2}$, as expected from Theorem~\ref{thm:inhPP}. Note that Theorem~\ref{thm:inhPP} provides a $(j \pi)^{-2}$ regime. The black line corresponds to the first bisector, so that the points align if the empirical convergence matches the theoretical regime. \label{Fig:PP_eigenval}}
\end{figure}

We first simulate $n=100$ Poisson processes using the \texttt{hawkesbow} package \cite{hawkesbow}. We consider different intensity functions: constant, $w(t) = 100$, linear, $w(t) = t$, $w(t) = \exp (-0.005 t)$, sinusoidal, $w(t) = (1+\sin(0.11 t))$. Then we compute the average eigenfunctions and eigenvalues over 50 replicates (Fig. \ref{Fig:PP_eigenfun} and \ref{Fig:PP_eigenval} respectively). For eigenvalues, the empirical asymptotic regime in $j$ conforms to the theoretical regime in $(j \pi)^{-2}$ (with a slight modification $(j \pi- \pi/2)^{-2}$ that better fits for small $j$'s, which matches the remarks following Theorem~\ref{thm:inhPP}). The form of the intensity does not impact much the convergence of eigenvalue estimates, even for large scales ($j$ small). As expected, eigenfunctions show increasing oscillations with $j$. However, the form of the intensity induces different oscillatory patterns in eigenfunctions  but the number of zero crossings is preserved with $j$ zeros on the interval $[0,1)$ for all functions.

	Then, to investigate the impact of dependencies in the data, we simulated $n=1,000$ Hawkes processes on $[0,100]$, with $w_0=100$, and self-exciting function $h(t) = \alpha e^{-\beta t}$. To explain briefly the simulation method, the package considers a reproduction function $h$ that is decomposed such as $h = m h^*$, with $m= \int_0^{\infty} h(t)dt$ and $h^*$ a true density function such that $\int_0^{\infty} h^*(t)dt = 1$. In our setting, this means that $h^*(t) = \beta \exp(-\beta t)$, and that $h(t) =\|h\|_1 \times h^*(t)$. Since Theorem~\ref{Hawkes-ed-sol} exhibits the central role of $\|h\|_1$ and $\alpha$, in our simulation we consider $\alpha \in \{0.01,0.1,1,10\}$ and $\|h\|_1 \in \{0.1,0.3,0.5,0.7,0.9\}$. Data were finally rescaled in $[0,1]$. This setting allows us to explore different strengths and ranges of dependency.  
	
\begin{figure}
	\begin{center}
		\includegraphics[scale=0.5]{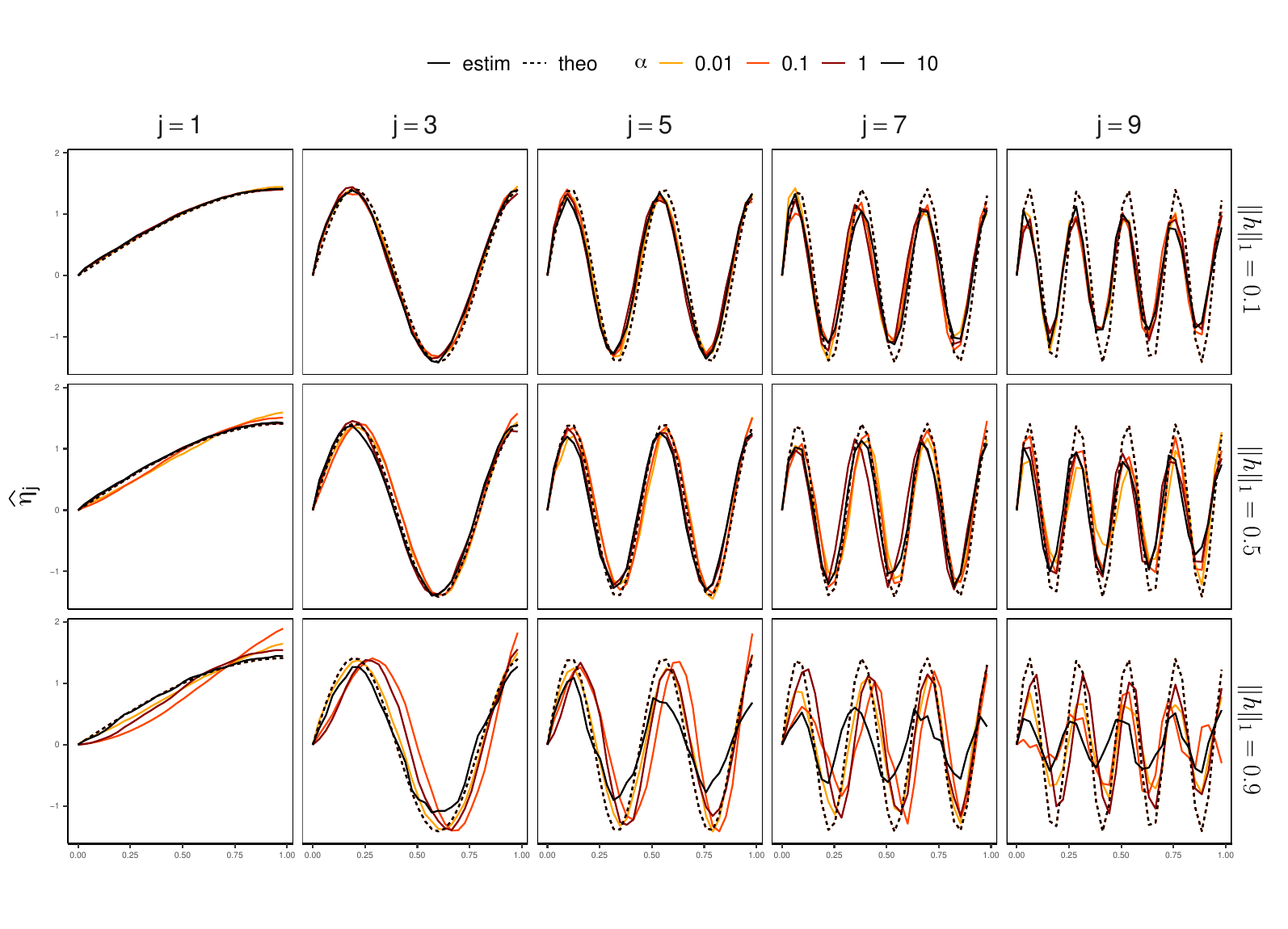}
	\end{center}
	\caption{Average eigenfunctions for Hawkes Processes with different transfert functions over 50 replicates. Dotted lines: asymptotic eigenfunctions $t\longmapsto\sqrt{2}\sin(\pi(2j-1)t/2)$ (see Theorem~\ref{Hawkes-ed-sol}). \label{Fig:H_eigenfun_alpha}}
\end{figure}

\begin{figure}
	\begin{center}
		\includegraphics[scale=0.8]{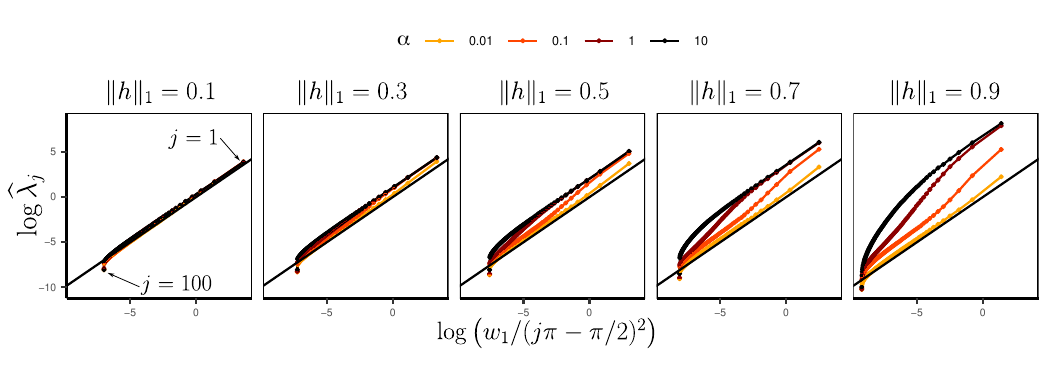}
	\end{center}
	\caption{Average eigenvalues (log-scale) for Hawkes Processes over 50 replicated. Each dot corresponds to a value of $j \in \{1,\hdots, 50\}$. The empirical average is plotted vs the expected theoretical asymptotic regime of eigenvalues in $w_1/(j \pi - \pi/2)^{2}$, as expected from Theorem~\ref{Hawkes-ed-sol}. The black line corresponds to the first bisector, so that the points align if the empirical convergence matches the theoretical regime. \label{Fig:H_eigenval_alpha}}
\end{figure}

 The eigenfunctions show a periodic pattern with zero crossings depending on $j$ (Fig.~\ref{Fig:H_eigenfun_alpha}). Interestingly, the number of zero crossings does not depend on the parameters $(\alpha,\|h\|_1)$ and agrees with the values predicted by the theory for all studied values of $j$ (note that for visualization purposes, eigenfunctions are studied for a small range of values of $j$). Fig.~\ref{Fig:H_eigenval_alpha} shows that convergence of eigenvalues holds for all parameters $(\alpha,\|h\|_1)$ but the convergence is slower when $\alpha$ or $\|h\|_1$ is large. Observe that in our setting, Assumption \eqref{cond-rho} is satisfied only in few cases: if $\alpha=0.01$, or if $\alpha=0.1$ and $\|h\|_1\in\{0.1,0.3\}$, so we conjecture that Assumption \eqref{cond-rho} may be a technical artefact and Theorem~\ref{Hawkes-ed-sol} may hold in more general situations.

From this study, we can conclude that many of the theoretical conclusions of Section~\ref{sec:illust} which are valid when $j$ is large seem to extend for small values of $j$ for Poisson and Hawkes processes (for the latter at least when $\|h\|_1$ and $\alpha$ are not too large). 
Note that the eigenfunctions for the Poisson and Hawkes processes have a very similar oscillating behavior, and this at each scale. 

\newpage
\section{Applications}\label{sec:applications}
\subsection{PCA on point processes reveals multiscale variability in earthquake occurrence data}
	
\begin{figure}
	\begin{center}
		\includegraphics[scale=0.6]{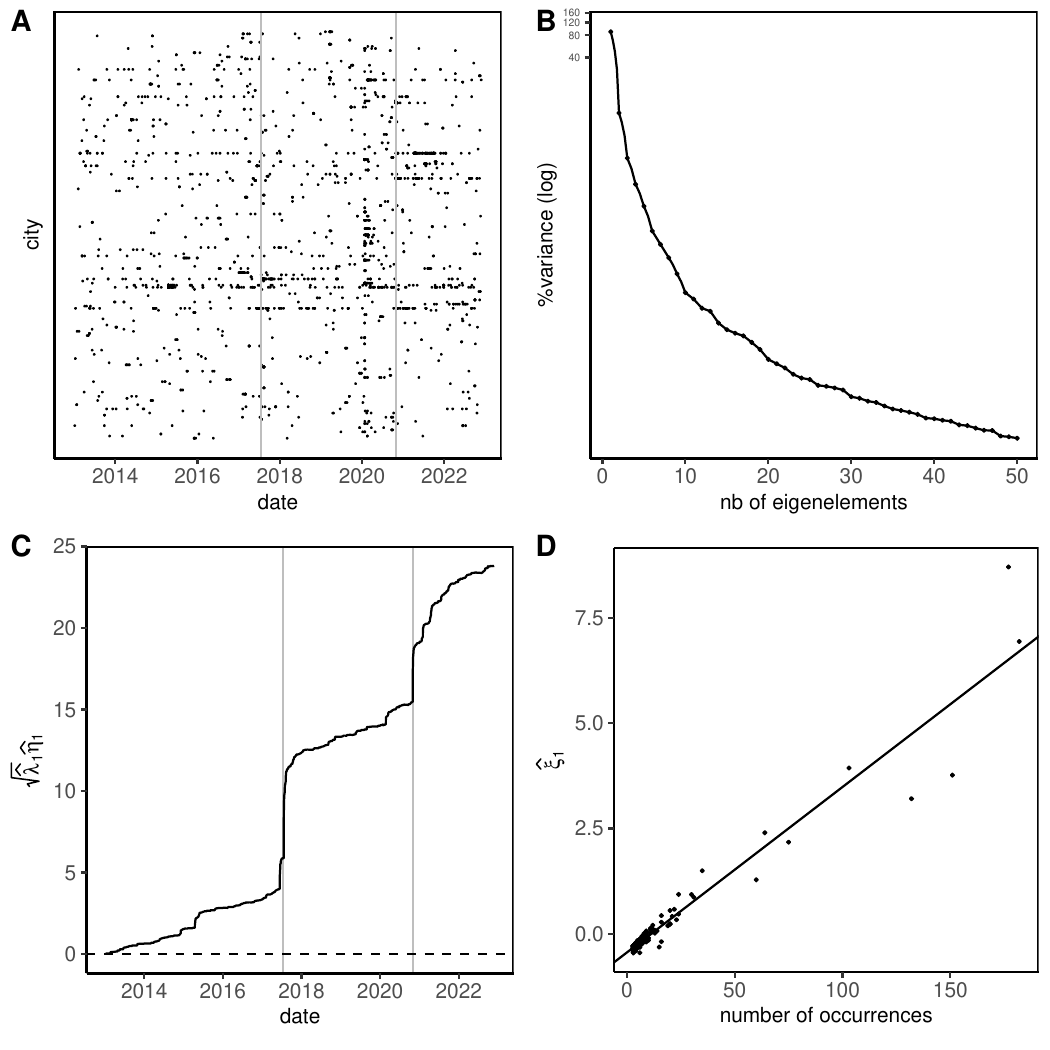}
	\end{center}
	\caption{A: Raster plot of earthquakes during the period 2013-2023. Each line corresponds to a city and each dot to an earthquake occurrence. Breakpoint dates (grey vertical lines) correspond to 2017.07.16 and 2020.11.01 B: Percentage of variance according to the number of eigenelements (log scale). C: First rescaled eigenfunction $\sqrt{\widehat{\lambda}_1}\widehat{\eta}_1$ according to the date of earthquakes (left). Breakpoint dates (grey vertical lines) correspond to 2017.07.16 and 2020.11.01. D: Plot of PCA scores for the first axis $\widehat\xi_{i1}$ according to the number of occurrences $W([0,t])$. \label{fig:Earthquake_axis1}}
\end{figure}

We begin by illustrating our method through the analysis of data obtained from the Kandilli Observatory and Earthquakes Research Institute at Bo\u{g}azi\c{c}i University\footnote{\url{http://www.koeri.boun.edu.tr/sismo/2/earthquake-catalog/}}. The dataset comprises earthquake occurrences in Turkey and neighboring regions of Greece, recorded between January 2013 and January 2023, spanning 1181 cities (Fig. \ref{fig:Earthquake_axis1}-A). In recent years, the Gulf of G\"{o}kova in Southwest Turkey has witnessed two significant seismic events: the Bodrum earthquake on July 16, 2017 (magnitude 6.6, also felt on the Greek island of Kos) \cite{Karasozen_etal_2018}, and the Aegean Sea earthquake on November~1, 2020, with a moment magnitude of 7.0 (the highest magnitude observed during the period).  In the following, we illustrate how our framework can be employed to investigate the fine-scale dynamics of earthquake occurrences in the region, focusing on the 195 cities that experienced more than two earthquakes over the specified period. Finally, note that our framework assumes the observed processes (earthquake occurrences) are independent from one city to another, even though spatial dependencies are likely to exist. The spatio-temporal modeling of such data lies beyond the scope of the present work.

We start by inspecting the significance of the first eigenfunction \(\widehat{\eta}_1\) (Fig. \ref{fig:Earthquake_axis1}-C) and the scores on this axis (Fig. \ref{fig:Earthquake_axis1}-D). We see that \(\widehat{\eta}_1\) corresponds to the temporal accumulation of the number of earthquakes across cities and over time, and that the scores \((\widehat{\xi}_{i1})_i\) directly correspond to the total number of earthquakes in each city over the specified period. Since this axis carries $\widehat{\lambda}_1/\sum_j \widehat{\lambda}_j = 89\%$ of the variability, we conclude that the first source of variance in the data is associated with the deviation of cities in the accumulation of earthquakes over time compared to the average temporal pattern of earthquakes in the region. We identify Akdeniz as an outlier with an unusually high number of earthquakes compared to the regional average. This city will be excluded from subsequent analyses that investigate the finer-scale dynamics of earthquake activity in the area.

\begin{figure}
	\centering
	\includegraphics[scale=0.8]{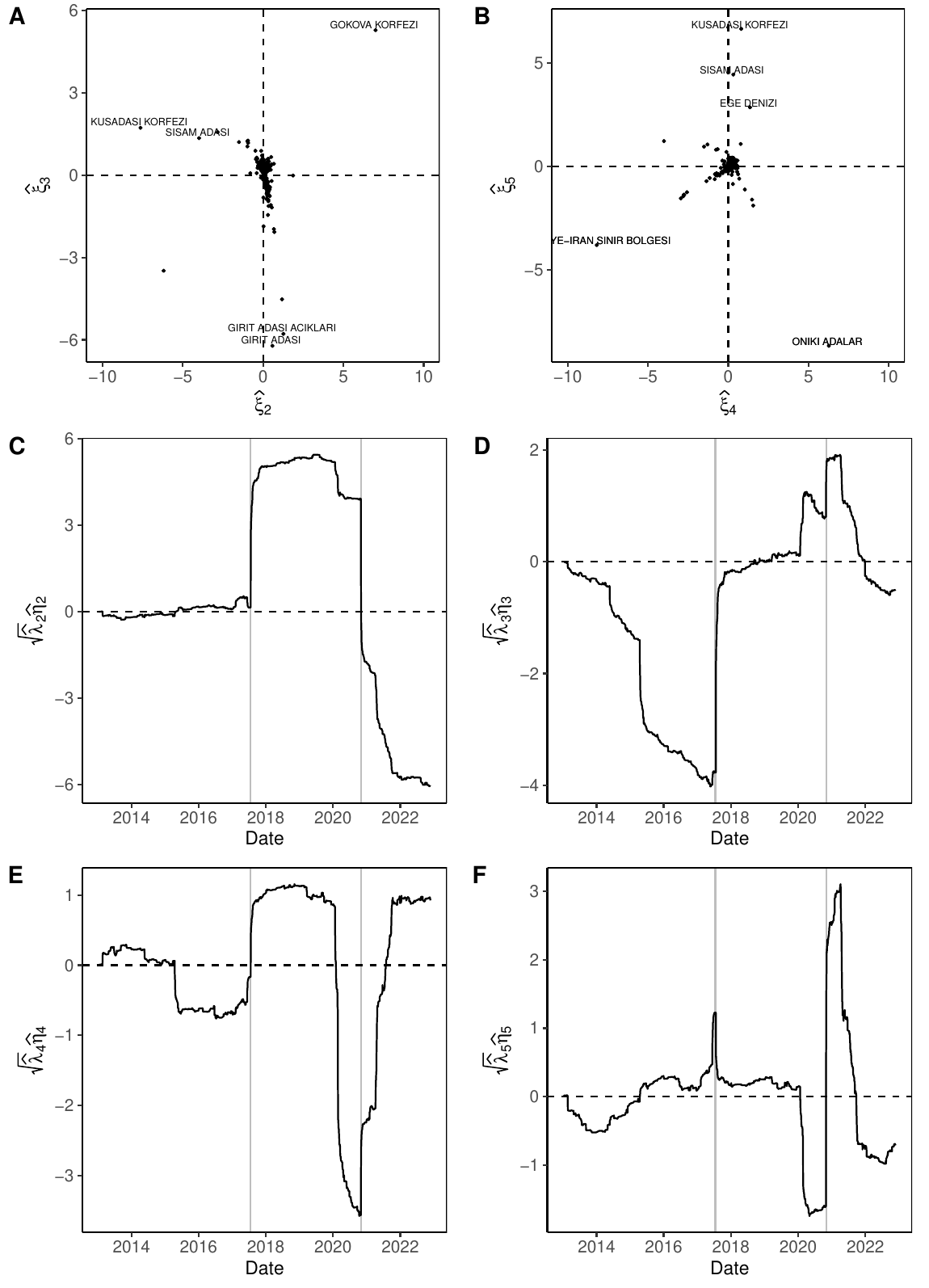}
	\caption{A: Estimated scores of cities $\widehat{\xi}_{i2}$ vs. $\widehat{\xi}_{i3}$. B: Estimated scores of cities $\widehat{\xi}_{i4}$ vs. $\widehat{\xi}_{i5}$. C-D-E-F: Estimated rescaled eigenfunctions $\sqrt{\widehat{\lambda}_j}\widehat{\eta}_j$ over time for $j=2,3,4,5$.\label{fig:Earthquake_scores_etas}}
\end{figure}

In order to interpret the subsequent axes, we recall that our method summarizes the dynamics of earthquake occurrences through simple functions, as expressed in Equation~\eqref{representation}. Consequently, we focus on the variations of the estimated eigenfunctions $(\widehat{\eta}_j)_j$, and represent the positions of cities according to their scores $(\widehat{\xi}_j)_j$ on successive axes as in any PCA analysis (Fig. \ref{fig:Earthquake_scores_etas}). These representations allow us to identify cities (like Gokova Korfezi, Onika Adalar) that have atypical accumulations of earthquakes over the period (these cities correspond to observations with atypical scores on the axes).

Interestingly, the second axis, \(\widehat\eta_2\) (Fig. \ref{fig:Earthquake_scores_etas}-C) reveals a distinct change in seismic activity between the two main earthquakes (2017 and 2020). Cities exhibiting positive scores on this axis (Fig. \ref{fig:Earthquake_scores_etas}-A) indicate a global increase in earthquake rates between July 2017 and November 2020, followed by a decrease below the average regional rate after November 2020. Conversely, cities with negative scores on \(\widehat{\eta}_3\) (Fig. \ref{fig:Earthquake_scores_etas}-A) show an accumulation of earthquakes before July 2017. Overall, our method provides a highly accurate description of the variability in earthquake occurrences among different cities. It offers a means to represent and position cities relative to each other based on their earthquake dynamics variability. Moreover, our estimation framework, relying on occurrence data without smoothing, captures sharp and fine-scale variations in this dynamics, as shown by the distinct peaks in $\widehat{\eta}_4$ and $\widehat{\eta}_5$ (Fig. \ref{fig:Earthquake_scores_etas}-E-~F).

\subsection{PCA on point processes to investigate single-cell epigenomic variations}

The field of single-cell biology has emerged as a significant framework for generating detailed molecular portraits of cell populations, providing valuable insights into the complexity of living tissues \cite{zheng_massively_2017}. While gene expression analysis has traditionally been instrumental in understanding fundamental biological processes, epigenomic-based gene regulation is receiving increasing attention, especially when mediated by chromatin modifications measured by sequencing-based chromatin immuno-precipitation assays (chIP-Seq). Briefly, chromatin comprises the fiber that constitutes chromosomes. ChIP-Seq assays enable the mapping of chromatin modifications, generating data in the form of 1D-spatial coordinates. However, there is currently no dedicated method for analyzing single-cell epigenomic profiles that adequately accounts for both the spatial nature of the signal and the single-cell resolution. 

A recent study \cite{marsolier_h3k27me3_2022} has shown that a particular histone modification (H3K27me3) is involved in the emergence of drug persistence in breast cancer cells. Drug persistence occurs when only a subset of cells, known as persister cells, survives the initial drug treatment, thereby creating a reservoir of cells from which resistant cells will emerge. Persister cells exhibit dynamic changes in H3K27me3 modifications at the single-cell level, but the variability of this pattern makes it difficult to relate to the emergence of resistance within tumor cells. Previous analyses have suggested that a pool of untreated cells could contribute to the persister cell population later upon exposure to chemotherapy. However, the lack of an appropriate methodological framework to handle the spatial characteristics of the data has made the clear identification of this pool of precursor cells difficult.

\begin{figure}
	\centering
	\includegraphics[scale=0.55]{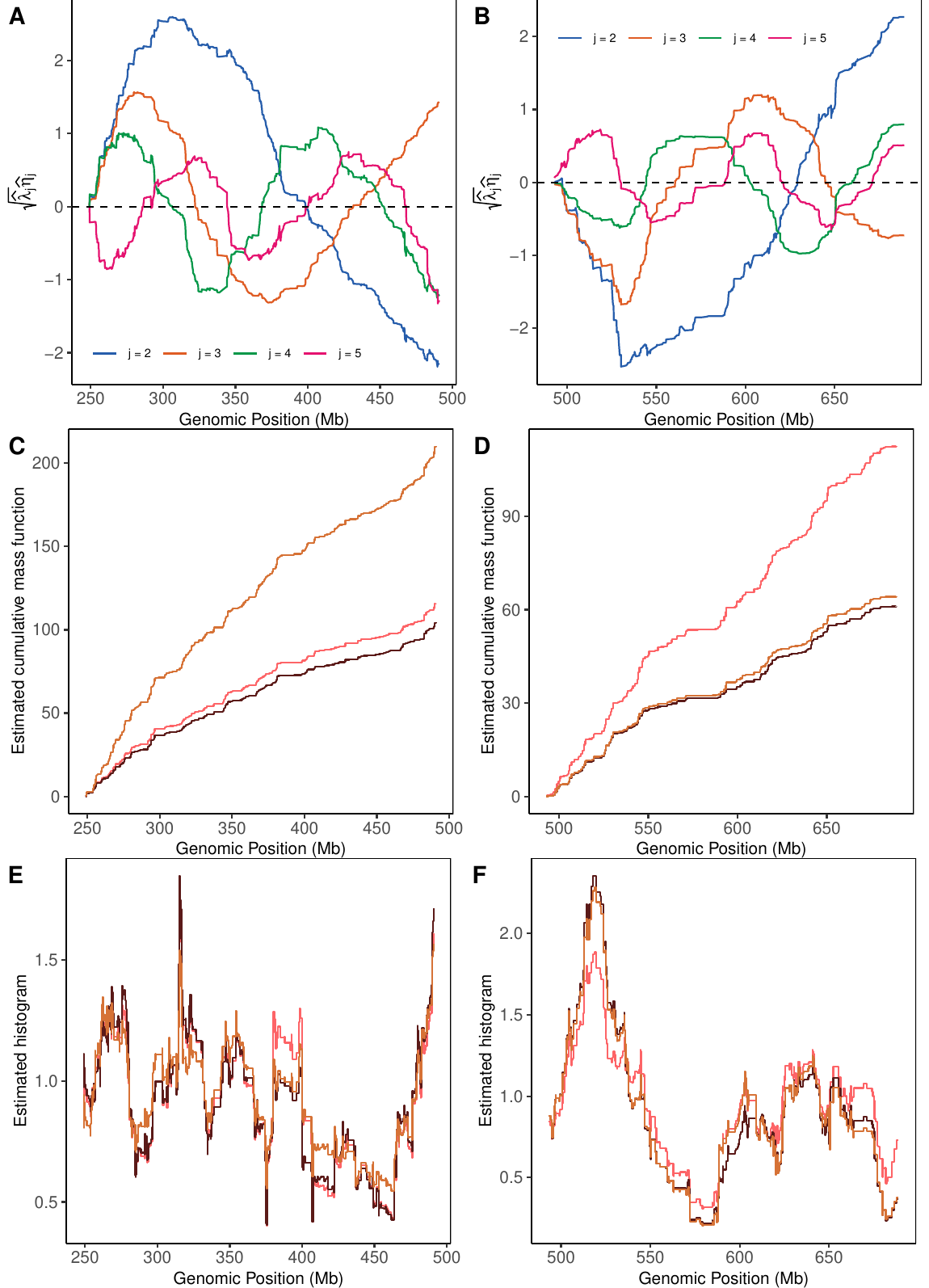}
	\caption{A-B: Eigenfunctions weighted by their respective eigenvalues $\sqrt{\widehat{\lambda}_j} \widehat{\eta}_j$ according to the position along the genome. Chromosome 2 (left) and Chromosome 3 (right). C-D: Estimated cluster-wise cumulative mass functions $
\widehat{\Pi}_k$. E-F: Estimated cluster-wise histograms $\widehat{w}_k$}. \label{fig:sc-chIPSeq}
\end{figure}

We propose using our PCA framework to explore the variability of H3K27me3 binding along the chromosomes 2 and 3 genomes among untreated cancer cells. To proceed, we performed PCA on $n=4722$ cells, with each cell described by the 1D spatial positions of H3K27me3 modifications along chromosomes. Eigenfunctions estimated by this approach allow us to visualize the basic components of spatial variability in H3K27me3 occurrences among untreated cancer cells (Fig. \ref{fig:sc-chIPSeq} A-B). Our PCA framework then provides principal scores $\widehat{\xi}_{ij}$, on which we conducted $k$-means clustering. Using $J=50$ principal components and selecting the number of clusters based on a silhouette criterion, we identified clusters of cells with similar profiles of histone modification accumulation across the genome. Since $k$-means clustering assigns label variables to each observed process ($\widehat{Z}_{ik}=1$ if cell $i$ is in cluster $k$, and $\widehat{n}_k=\sum_{i}\widehat{Z}_{ik}$ the size of cluster $k$), we represent the estimated cluster-wise cumulative mass function ($\widehat{\Pi}_k([0,t])$, Fig. \ref{fig:sc-chIPSeq} C-D) and associated histograms ($\widehat{w}_k(t)$, Fig. \ref{fig:sc-chIPSeq} E-F): 
$$
\widehat{\Pi}_k([0,t]) = \frac{1}{\widehat{n}_k} \sum_{i = 1}^{n} \widehat{Z}_{ik} \Pi_i([0,t]), \quad \widehat{w}_k(t) = \frac{1}{\widehat{n}_k} \sum_{i = 1}^{n}  \widehat{Z}_{ik}\sum_{T \in N_i} \mathbf{1}_{\{ T=t \}}.
$$
This analysis reveals the spatial variability of histone modifications among cancer cells, offering a new perspective on epigenomic heterogeneity by revealing three sub-populations within untreated cancer cells, characterized by zones along chromosomes that concentrate higher or lower numbers of occurrences. Further research is needed to fully understand the biological implications of these variabilities and their role in forming a pool of persister cells. Nonetheless, our PCA for point processes framework marks an important step towards investigating spatial epigenomic variations among cell populations.

\subsection{PCA on point process for multi-electrodes spike data}

Neuroscience has witnessed high-paced developments in recent years, thanks to the convergence of experimental and computational approaches. Signal processing and machine learning, in general, are deeply connected to the analysis of neuroscience data, often in the form of spike trains that correspond to neuronal activities. These spike trains are typically modeled as point processes, and technological advancements now allow for the recording of such spikes in populations of neurons, paving the way for the analysis framework of multiple neural spike train data \cite{Cunningham2014}. Among many methodological challenges, dimension reduction emerges as a key challenge for visualizing, exploring, and understanding the structure and variability of population activity. The goal of dimensionality reduction is to characterize how the firing rates of different neurons co-vary, while discarding the spiking variability as noise. We present how PCA on point processes can be utilized in this context. In particular, the principal measure we infer can be interpreted as a latent common point process from which individual action potentials are generated in a stochastic manner. Hence, our representation with principal measures enables us to define a latent space that represents shared activity patterns prominent in the population response.

We consider data from an optogenetic therapy experiment performed on non-human primates retina. The aim of this experiment is to restore light sensitivity in the residual retinal tissue after photoreceptor generation like macular degeneration. These diseases affect photoreceptor cells, but the retinal ganglion cells (RGC) can still communicate with the brain via the optic nerve. Then the data consists in the recording the neuronal activity of retinal fragments placed on Multi-Electrodes-Array (MEA, $n=256$ electrodes) grid. These 256 electrodes form a 16 x 16 grid onto a device that is positioned on the retina of macaques. The activity of each electrode can be related to a particular retina region, so that the data are typically spatio-temporal, but we only consider its temporal components. Different visual stimuli are applied to retinal fragments and the activity of electrodes is recorded. 

\begin{figure}
	\centering
	\includegraphics[scale=0.55]{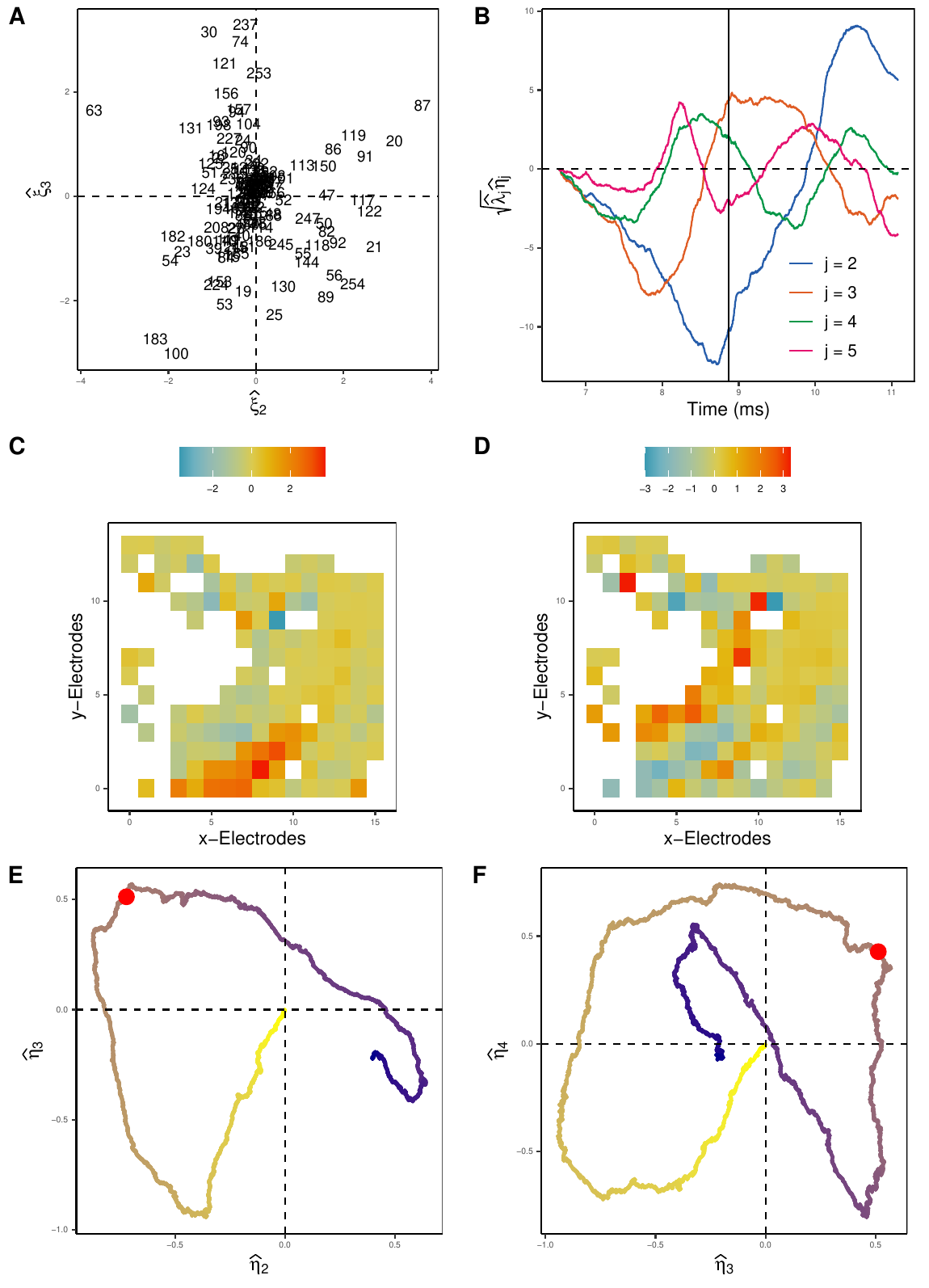}
	\caption{A: scores $\widehat{\xi}_2$ vs $\widehat{\xi}_3$ for electrodes $n=1, \hdots, 256$. B: Eigenfunctions weighted by their respective eigenvalues $\sqrt{\widehat{\lambda}}_j \widehat{\eta}_j$ according to time (in ms). The vertical line corresponds to the stimulus at $t=8.86$ms. C-D: visualisation of scores $\widehat{\xi}_2$ (C) and  $\widehat{\xi}_3$ (D) according to the position of each electrode on the retinal device. Red: high value, blue: low value. E-F: representation of the data in the latent space (E, $\widehat{\eta_3}$ vs $\widehat{\eta_2}$, and F, $\widehat{\eta_4}$ vs $\widehat{\eta_3}$). Time is represented by colors (yellow: begining of the assay, purple: end of the assay). Red dot: stimulus time.\label{fig:neuro}}
\end{figure}

Thanks to our framework, we can easily visualize the variability in the response of each electrode (Fig. \ref{fig:neuro} A), as well as the temporal dynamics that explain this variability in spike train occurrences (Fig. \ref{fig:neuro} B). These principal components can be interpreted as latent spike trains that compose the dynamics of the population. These latent processes could be of interest for understanding the neuronal coordination underlying the population response. Additionally, the score of each electrode can be represented with respect to the position of the electrodes on the device (Fig. \ref{fig:neuro} C-D) to visualize variations in each principal component across the device. Identifying regions on the device that show high or low scores could be related to brain regions whose responses are governed by the corresponding latent spiking process. Finally, this dynamical system can also be summarized by representing trajectories in the principal measure space (Fig. \ref{fig:neuro} E-F). In this representation, each time point $t$ corresponds to a single point in the latent firing rate space $(\widehat{\eta}_2(t), \widehat{\eta}_3(t))$. This representation allows us to follow the temporal dynamics of the observed neuron population. Hence, we expect that PCA on point processes can become a natural framework to decompose these complex phenomena of population-based spike train analysis by providing simple decompositions that help unravel the diversity and variability of neuronal population spiking activities. In this application, the assumption of independence between neurons may seem as a strong one. Exploring the spatial dependency between neurons could be an interesting research direction to enhance the method.

\section{Perspectives}
In this work, we propose a novel mathematical and statistical framework for dimension reduction in replicated one-dimensional point processes. Our approach is based on functional PCA for random measures, for which we derive a Karhunen-Loève expansion that enables us to define the concept of  principal measure. The key tool is the associated covariance operator whose eigenfunctions are the cumulative mass functions of the principal measures. We provide a detailed characterization of these eigenelements in simple point process models such as the Poisson and Hawkes models. Assuming that the point processes are independent and identically distributed, we further develop an easily implementable, fully data-driven estimation procedure for these eigenelements that requires no smoothing. The resulting estimators achieve the parametric rate of convergence. These rates may, however, fail to hold when dependence is present across replicates (e.g., across spatial locations or neurons). Extending the methodology to dependent settings would be an interesting direction for future research. Our contribution also includes various applications across different fields, as well as a software package that implements our methodology.  \\
This work opens the door to many potential extensions, with a particularly interesting direction being the adaptation of our framework to multivariate processes. For example, the observed processes could be spatial point patterns, where \((N_1, \hdots, N_n)\) are observed on \([0,1]^d\). In this case, our framework could be applied to analyze image data, for instance, with \(d = 2\) or \(d = 3\). Another scenario involves each process \(N_i\) being observed on \([0,1]\), where the observed data consist of multiple dependent processes, such that \(N_i = (N_i^1, \hdots, N_i^d)\). Our methodology could then be extended to multivariate functional PCA, which seeks a joint representation of the \(d\) dependent processes \cite{Happ:Greven}. In neuroscience, for example, repeated assays over time often involve multiple dependent neurons. Multivariate functional PCA for multivariate point processes could help uncover latent components that explain the variability in spiking events across time and neurons. Similarly, in single-cell genomics, new technologies allow for the simultaneous recording of multiple dependent chromatin modifications within a single assay at the single-cell level \cite{Meers2023}. Extending our framework to multivariate processes requires a separate and dedicated study, in particular because our current methodology relies heavily on the concept of the cumulative mass function \(F_{\mu}\), which is defined on the one-dimensional interval \([0,1]\) and takes values in \(\mathbb{R}\). Extending our approach will thus necessitate the development of specific mathematical tools tailored to higher-dimensional settings.
 
\begin{acks}[Acknowledgments]
	The authors would like to thank Emmanuel Tr\'elat for valuable suggestions concerning solutions of Sturm-Liouville theory, Celine Vallot for discussions on the single-cell epigenomic data, and Gregory Gauvain, Fabrice Arcizet and Pierre Pouget for sharing the neuroscience data. The authors would like to warmly thank the Associate Editor and the anonymous referees for very valuable comments and suggestions.
\end{acks}

\begin{funding}
	The research was supported by a grant from the Agence Nationale de la Recherche ANR-18-CE45-0023 SingleStatOmics, by the projects AI4scMed, France 2030 ANR-22-PESN-0002.
\end{funding}

\bibliographystyle{imsart-nameyear}

\section{Proofs of Section~\ref{sec:KLMercer}}
\subsection{Proof of Proposition~\ref{prop:muj}}
Let $j\geq 1$. We observe that, since $\lambda_j>0$, 
\begin{equation}\label{eq:kernel_eig}
\eta_j(t)=\lambda_j^{-1}\int_0^1K_\Delta(t,s)\eta_j(s)ds, \qquad t\in[0,1]. 
\end{equation}
We prove that $\eta_j$ is of finite variation on $[0,1]$.
We recall that a function $f$ is a function of finite variation on an interval $[a,b]$ if, for any grid $(t_0,\hdots,t_p)$ of $[0,1]$, the quantity 
\[
V_f := \sup_{a\leq t_0<t_1<\hdots<t_p\leq b}\left\{\sum_{\ell=1}^p|f(t_\ell)-f(t_{\ell-1})|\right\}
\]
is finite. For any grid $(t_0,\hdots,t_p)$ of $[0,1]$, we have:
\begin{eqnarray*}
\sum_{\ell=1}^p|\eta_j(t_\ell)-\eta_j(t_{\ell-1})|&=&\lambda_j^{-1}\sum_{\ell=1}^p\left|\int_0^1(K_\Delta(s,t_\ell)-K_\Delta(s,t_{\ell-1}))\eta_j(s)ds\right|\\
&\hspace{-4cm}=&\hspace{-2cm}\lambda_j^{-1}\sum_{\ell=1}^p\left|\int_0^1\left(\mathbb E[F_\Delta(s)F_\Delta(t_\ell)]-\mathbb E[F_\Delta(s)F_\Delta(t_{\ell-1})]\right)\eta_j(s)ds\right|\\
&\hspace{-4cm}=&\hspace{-2cm}\lambda_j^{-1}\sum_{\ell=1}^p\left|\int_0^1\mathbb E[F_\Delta(s)(F_\Delta(t_{\ell})-F_\Delta(t_{\ell-1}))]\eta_j(s)ds\right|\\
&\hspace{-4cm}\leq &\hspace{-2cm}\lambda_j^{-1}\int_0^1\mathbb E\left[\left|F_\Delta(s)\right|\sum_{\ell=1}^p\left|F_\Delta(t_{\ell})-F_\Delta(t_{\ell-1})\right|\right]\left|\eta_j(s)\right|ds\\
&\hspace{-4cm}\leq &\hspace{-2cm}\lambda_j^{-1}\int_0^1\mathbb E\left[\left|F_\Delta(s)\right|\sum_{\ell=1}^p\left|\Delta(]t_{\ell-1},t_{\ell}])\right|\right]\left|\eta_j(s)\right|ds\\
&\hspace{-4cm}\leq &\hspace{-2cm}\lambda_j^{-1}\int_0^1\mathbb E\left[\left|F_\Delta(s)\right|\sum_{\ell=1}^p\left(\Pi(]t_{\ell-1},t_{\ell}])+W(]t_{\ell-1},t_{\ell}])\right)\right]\left|\eta_j(s)\right|ds\\
&\hspace{-4cm}\leq &\hspace{-2cm}\lambda_j^{-1}\int_0^1\mathbb E\left[\left|F_\Delta(s)\right|\left(\Pi([0,1])+W([0,1])\right)\right]\left|\eta_j(s)\right|ds,
\end{eqnarray*}
which is independent of the grid. We observe that
\begin{align*}
\mathbb E\left[\left|F_\Delta(s)\right|\left(\Pi([0,1])+W([0,1])\right)\right]&\leq \sqrt{\mathbb E\left[\left|F_\Delta(s)\right|^2\right]\times\mathbb E\left[\left(\Pi([0,1])+W([0,1])\right)^2\right]}\\
&\leq 2\sqrt{\mathbb E\left[\left|F_\Delta(s)\right|^2\right]\mathbb E\left[\Pi^2([0,1])\right]}
\end{align*}
and
\begin{align*}
\int_0^1\mathbb E\left[\left|F_\Delta(s)\right|\left(\Pi([0,1])+W([0,1])\right)\right]\left|\eta_j(s)\right|ds&\\
&\hspace{-3cm}\leq 2\|\eta_j\|_2\sqrt{\mathbb E\left[\Pi^2([0,1])\right]}\Bigg(\int_0^1\mathbb E\left[\left|F_\Delta(s)\right|^2\right]ds\Bigg)^{1/2}\\
&\hspace{-3cm}\leq 4\|\eta_j\|_2\mathbb E\left[\Pi^2([0,1])\right]<\infty.
\end{align*}
Then, $\eta_j$ is a function of finite variation. Since $K_\Delta$ is continuous, $\eta_j$ is also continuous (and then right-continuous). We also have that $\eta_j(0)=0$. Then, applying Proposition~4.4.3 of \cite{Cohn93}, we obtain that there exists a signed measure $\mu_j$ such that
\[
\eta_j(t) = \mu_j([0,t]),\quad t\in [0,1].
\]
From Fubini's Theorem we have that for any function $\varphi\in\mathcal H_0^1$ 
\begin{align}\label{phimu}
\int_0^1 \varphi'(t)\eta_j(t)dt&=\int_0^1\varphi'(t)\Bigg(\int_0^td\mu_j(s)\Bigg)dt=\int_0^1\Bigg(\int_s^1\varphi'(t)dt\Bigg)d\mu_j(s)\nonumber\\
&=-\int_0^1 \varphi(s)d\mu_j(s),
\end{align}
since $\varphi(1)=0$.
\subsection{Proof of Theorem~\ref{thm:KL}}
By the Karhunen-Lo\`eve theorem (see e.g. \cite{hsing2015theoretical}, Theorem 7.3.5) applied to the stochastic process $F_\Pi$,
\begin{equation}\label{eq:Mercer2}
\lim_{J\to+\infty}\sup_{t\in [0,1]}\mathbb E\left[\left(F_\Pi(t)-F_W(t)-\sum_{j=1}^J\sqrt{\lambda_j}\xi_j\eta_j(t)\right)^2\right]=0.
\end{equation}
Furthermore, similarly to \eqref{phimu}, we have, for all $\varphi\in\mathcal H_0^1$, 
$$\langle \varphi,\Pi-W\rangle=-\langle F_\Pi-F_W,\varphi'\rangle.$$
Then, using Equation~\eqref{phimu*} of the main file,
\begin{align*}
\mathbb E\Bigg[\Big\|\Pi-W-\sum_{j=1}^J\sqrt\lambda_j\xi_j\mu_j\Big\|_{\mathcal H^{-1}}^2\Bigg]&\leq\mathbb E\Bigg[\sup_{\varphi\in\mathcal H_0^1, \|\varphi'\|\leq 1}\Big(\langle \varphi,\Pi-W\rangle-\sum_{j=1}^J\sqrt{\lambda_j}\xi_j\langle \varphi,\mu_j\rangle\Big)^2\Bigg]\\
&\leq\mathbb E\Bigg[\sup_{\varphi\in\mathcal H_0^1, \|\varphi'\|\leq 1}\Big\langle F_\Pi-F_W-\sum_{j=1}^J\sqrt{\lambda_j}\xi_j \eta_j,\varphi'\Big\rangle^2\Bigg]\\
&\leq \mathbb E\Bigg[\sup_{\varphi\in\mathcal H_0^1, \|\varphi'\|\leq 1}\|\varphi'\|^2\Big\| F_\Pi-F_W-\sum_{j=1}^J\sqrt{\lambda_j}\xi_j \eta_j\Big\|^2\Bigg]\\
&\leq \mathbb E\Bigg[\int_0^1\Big( F_\Pi(t)-F_W(t)-\sum_{j=1}^J\sqrt{\lambda_j}\xi_j \eta_j(t)\Big)^2dt\Bigg]
\end{align*}
and the result comes from  Equation~\eqref{eq:Mercer2}.

\subsection{Proof of Theorem~\ref{thm:Mercer}}


A direct consequence of Mercer's theorem is that
\[
\left\langle K_\Delta, \frac{\partial^2 \varphi}{\partial t\partial s}\right\rangle = \left\langle \sum_{j\geq 1}\lambda_j\eta_j\otimes\eta_j, \frac{\partial^2 \varphi}{\partial t\partial s}\right\rangle =\sum_{j\geq 1}\lambda_j \left\langle \eta_j\otimes\eta_j, \frac{\partial^2 \varphi}{\partial t\partial s}\right\rangle=\sum_{j\geq 1}\lambda_j \left\langle \varphi, \mu_j\otimes\mu_j\right\rangle. 
\]
Indeed,  Fubini's theorem and the definition of $\mu_j$ implies 
\begin{eqnarray*}
\left\langle \eta_j\otimes\eta_j, \frac{\partial^2 \varphi}{\partial t\partial s}\right\rangle &=&\int_{[0,1]^2}\eta_j(s)\eta_j(t)\frac{\partial^2 \varphi}{\partial t\partial s}(s,t)dsdt\\
&=&\int_0^1\eta_j(s)\left(\int_0^1\eta_j(t)\frac{\partial^2 \varphi}{\partial t\partial s}(s,t)dt\right)ds\\
&=&\int_0^1\eta_j(s)\left\langle\eta_j,\frac{\partial^2 \varphi}{\partial t\partial s}(s,\cdot)\right\rangle ds\\
&=&-\int_0^1\eta_j(s)\left\langle\frac{\partial \varphi}{\partial s}(s,\cdot),\mu_j\right\rangle ds\\
&=&-\int_0^1\eta_j(s)\int_0^1\frac{\partial \varphi}{\partial s}(s,t)d\mu_j(t) ds\\
&=&\langle\varphi,\mu_j\otimes\mu_j\rangle,
\end{eqnarray*}
using Fubini's theorem and the definition of $\mu_j$ a second time. 
Hence we can write 
\[
\left\langle\frac{\partial^2 K_\Delta}{\partial t\partial s},\varphi\right\rangle = \sum_{j\geq 1}\lambda_j\langle\varphi,\mu_j\otimes\mu_j\rangle. 
\]
Now remark also that,
\[
C_\Delta = \frac{\partial^2 K_\Delta}{\partial t\partial s}. 
\] 
Indeed, for a function $\varphi\in\mathcal H_0^1$, we get 
\begin{eqnarray*}
\int_0^1F_\Pi(t)\varphi'(t)dt &=& \int_0^1\sum_{T\in N}\mathbf 1_{T\leq t}\varphi'(t)dt=\sum_{T\in N}\int_T^1\varphi'(t)dt\\
&=&\sum_{T\in N}(\varphi(1)-\varphi(T))=-\langle\Pi,\varphi\rangle. 
\end{eqnarray*}
In particular 
\[
\mathbb E\left[\int_0^1F_\Pi(t)\varphi'(t)dt\right] = -\mathbb E[\langle\Pi,\varphi\rangle] = -\langle W,\varphi\rangle. 
\]
This implies
\begin{eqnarray*}
\int_{[0,1]^2}K_\Delta(s,t) \frac{\partial^2 \varphi}{\partial t\partial s}(s,t)dsdt &=& \int_{[0,1]^2}\mathbb E[F_\Delta(s)F_\Delta(t) ]\frac{\partial^2 \varphi}{\partial t\partial s}(s,t)dsdt\\
&=&\mathbb E\left[\int_{[0,1]^2}F_\Pi(s)F_\Pi(t)\frac{\partial^2 \varphi}{\partial t\partial s}(s,t)dsdt\right]\\
&&-\int_{[0,1]^2}F_W(s)F_W(t)\frac{\partial^2 \varphi}{\partial t\partial s}(s,t)dsdt\\
&=&\mathbb E\left[\int_0^1F_\Pi(s)\left(\int_0^1F_\Pi(t)\frac{\partial^2 \varphi}{\partial t\partial s}(s,t)dt\right)ds\right]\\
&&-\int_0^1F_W(s)\mathbb E\left[\int_0^1F_\Pi(t)\frac{\partial^2 \varphi}{\partial t\partial s}(s,t)dt\right]ds\\
&=&-\mathbb E\left[\int_0^1F_\Pi(s)\sum_{T\in N}\frac{\partial \varphi}{\partial s}(s,T)ds\right]\\
&&+\int_0^1F_W(s)\int_0^1\frac{\partial \varphi}{\partial s}(s,t)W(t)dtds\\
&=&\mathbb E\left[\sum_{T,T'\in N}\varphi(T',T)\right]-\langle \varphi,W\otimes W\rangle\\
&=&\langle C_\Delta,\varphi\rangle. 
\end{eqnarray*}
Then, for all $\varphi\in\mathcal H_0^2$, from Cauchy-Schwarz's inequality we deduce 
\[
\left|\left\langle\sum_{j=1}^J\lambda_j\mu_j\otimes\mu_j-C_\Delta,\varphi\right\rangle\right| = \left|\left\langle \sum_{j=1}^J\lambda_j\eta_j\otimes\eta_j-K_\Delta,\frac{\partial^2\varphi}{\partial s\partial t}\right\rangle \right|\leq \left\|K_\Delta-\sum_{j=1}^J\lambda_j\eta_j\otimes\eta_j\right\|\left\|\frac{\partial^2\varphi}{\partial s\partial t}\right\|,
\]
which implies
\[
\left\|\sum_{j=1}^J\lambda_j\mu_j\otimes\mu_j-C_\Delta\right\|_{\mathcal H^{-2}}\leq \left\|\sum_{j=1}^J\lambda_j\eta_j\otimes\eta_j-K_\Delta\right\|, 
\]
and Mercer theorem implies the expected result. 

\section{Proofs of Section~\ref{sec:illust}}
\subsection{Proof of Theorem~\ref{KL-SLP}}
First, assume that $F_{\mu_j}$ is an eigenfunction of the operator $\Gamma_\Delta$ associated to the eigenvalue $\lambda_j$. We have: 
\begin{eqnarray}\label{poisson-start}
\lambda_j F_{\mu_j}(t)&=&\int_0^1K_\Delta(t,s)F_{\mu_j}(s)ds,\quad t\in[0,1].
\end{eqnarray}
We have for any $(t,t')\in[0,1]^2$, since $\lambda_j>0$,
$$|F_{\mu_j}(t)-F_{\mu_j}(t')|\leq \|F_{\mu_j}\|_2\times \lambda_j^{-1}\Bigg(\int_0^1\big(K_\Delta(t,s)-K_\Delta(t',s)\big)^2ds\Bigg)^{1/2}.$$
And since $K_\Delta$ is continuous, $|F_{\mu_j}(t)-F_{\mu_j}(t')|$ goes to $0$ when $t'\to t$. So $F_{\mu_j}$ is continuous.
Equation \eqref{poisson-start} gives
\begin{eqnarray}\label{poisson-equa}
\lambda_j F_{\mu_j}(t) &=& \int_0^1\Big(\int_0^{\min\{s;t\}} w(u)du\Big)F_{\mu_j}(s)ds\nonumber\\
&=&\int_0^tF_{\mu_j}(s)\Big(\int_0^s w(u)du\Big)ds+\int_t^1F_{\mu_j}(s)\Big(\int_0^t w(u)du\Big)ds.
\end{eqnarray}
Now, consider
\[
H_j(t)=\int_t^1F_{\mu_j}(s)ds,\quad t\in [0,1].
\]
Since $F_{\mu_j}$ is continuous, $H_j$ is a differentiable function and
\[
H'_j(t)=-F_{\mu_j}(t).
\]
 Actually, previous computations show that $F_{\mu_j}$ is an infinitely-differentiable function and differentiating Equation~\eqref{poisson-equa} gives
\begin{eqnarray} \label{eq:PPKernel}
\lambda_jF_{\mu_j}'(t) &=& w(t)\int_t^1F_{\mu_j}(u)du
\end{eqnarray}
and then $H_j$ is solution of a second-order differential equation, 
\begin{eqnarray}\label{eq:equa_diff-PP}
-\lambda_jH_j''(t) &= &H_j(t) w(t),
\end{eqnarray}
with boundary conditions:
\[
H_j(1)=0,\quad H_j'(0)=-F_{\mu_j}(0)=0.
\]
Therefore, $(\lambda_j,H_j)$ is solution of the system~\eqref{SLP} of the main file.

\bigskip 

Conversely, assume that $(\lambda_j,H_j)$ is a solution of \eqref{SLP} of the main file.  Then, setting
$$F_{\mu_j}(s)=-H'_j(s),\quad s\in[0,1],$$
we obtain, by using integration by part, for any $t\in[0,1]$,
\begin{align*}
\int_0^1K_\Delta(t,s)F_{\mu_j}(s)ds&=\int_0^tF_{\mu_j}(s)\Big(\int_0^s w(u)du\Big)ds+\int_t^1F_{\mu_j}(s)\Big(\int_0^t w(u)du\Big)ds\\
&=-\int_0^tH'_j(s)\Big(\int_0^s w(u)du\Big)ds-\Big(\int_0^t w(u)du\Big)\int_t^1H'_j(s)ds\\
&=-H_j(t)\int_0^tw(u)du+\int_0^tH_j(s)w(s)ds-\Big(\int_0^t w(u)du\Big)\Big(H_j(1)-H_j(t)\Big)\\
&=\int_0^tH_j(s)w(s)ds,
\end{align*}
since $H_j(1)=0$. Then, since $H'_j(0)=0$, for any $t\in[0,1]$,
\begin{align*}
\int_0^1K_\Delta(t,s)F_{\mu_j}(s)ds&=-\lambda_j\int_0^tH''_j(s)ds\\
&=-\lambda_jH'_j(t)=\lambda_jF_{\mu_j}(t).
\end{align*}
Then, 
$(\lambda_j,F_{\mu_j})$ are eigenelements of the operator of kernel $K_\Delta$.

\subsection{Proofs of Theorem~\ref{Hawkes-ed} and of the result of Remark~\ref{eigenvalue=0Hawkes}}
In the sequel, we set
\[
\theta=\beta-\alpha,\quad C_1=\frac{\beta^3w_0}{(\beta-\alpha)^3},\quad C_2=-\frac{w_0\alpha\beta(2\beta-\alpha)}{2(\beta-\alpha)^4}.
\] 
We observe 
\begin{align*}
(\alpha-\beta)^{-2}(2C_2(\alpha-\beta)-C_1)&=(\alpha-\beta)^{-2}\Big(\frac{w_0\alpha\beta(2\beta-\alpha)}{(\beta-\alpha)^3}-\frac{\beta^3w_0}{(\beta-\alpha)^3}\Big)\\
&=w_0\beta(\beta-\alpha)^{-5}\big(\alpha(2\beta-\alpha)-\beta^2\big)\\
&=-\frac{w_0\beta}{(\beta-\alpha)^3}=-C_1\beta^{-2}
\end{align*}
and 
\begin{equation}\label{C1C2}
2C_2\theta +C_1=C_1\beta^{-2}\theta^2.
\end{equation}
Using Equation~\eqref{eq:K_Hawkes_exp} of the main file, we have, by symmetry of the bivariate kernel
\[
K_\Delta(s,t)=C_1\min\{s;t\}-C_2\Big(-1-e^{-\theta|t-s|}+e^{-\theta t}+e^{-\theta s}\Big),\quad s,t\in[0,1].
\]
Furthermore, using definitions of $C_1$ and $C_2$, $(\lambda_j,H_j)$ is solution of \eqref{SL-Hawkes} of the main file if and only if
we have:
\[
\left\{
\begin{array}{ll}
\displaystyle -\lambda_j H''_j(s)=C_1\beta^{-2}\theta^2H_j(s)-C_2\theta^2 \int_0^1 e^{-\theta|t-s|}H_j(t)dt,   &s\in(0,1),\\
H_j(1)  =0,\quad H'_j(0)=0&
\end{array}
\right.
\]
Then, let us express $\Gamma_\Delta F_{\mu_j}$. For any $s\in[0,1]$,
\begin{eqnarray}
\Gamma_\Delta F_{\mu_j}(s) &=& \int_0^1K_\Delta(s,t) F_{\mu_j}(t)dt\nonumber\\
&=&\int_0^1\left(C_1\min\{s,t\}-C_2\left(-1-e^{-\theta|t-s|}+e^{-\theta t}+e^{-\theta s}\right)\right) F_{\mu_j}(t)dt\nonumber\\
&=&C_1\int_0^s t F_{\mu_j}(t)dt+C_1s\int_s^1  F_{\mu_j}(t)dt+C_2(1-e^{-\theta s})\int_0^1 F_{\mu_j}(t)dt \nonumber\\
&&-C_2\int_0^1e^{-\theta t} F_{\mu_j}(t)dt+C_2e^{\theta s}\int_s^1e^{-\theta t} F_{\mu_j}(t)dt+C_2e^{-\theta s}\int_0^se^{\theta t} F_{\mu_j}(t)dt.\nonumber
\end{eqnarray}
We have
\[
H_j(s)=\int_s^1 F_{\mu_j}(t)dt,\quad s\in [0,1].
\]
Since $K_\Delta$ is continuous, the eigenfunctions are continuous as well (eigenvalues are positive). 
Since $ F_{\mu_j}$ is continuous, we have
\[
 F_{\mu_j}(s)=-H'_j(s),\quad s\in [0,1]
\]
with boundary conditions
\[
H_j(1)=0,\quad H'_j(0)=0.
\]
Then, we can write
\[
\Gamma_\Delta F_{\mu_j}(s) =(I)+(II)-(III)+(IV)+(V)
\]
with
\begin{align*}
(I)&=C_1\int_0^s t F_{\mu_j}(t)dt+C_1s\int_s^1  F_{\mu_j}(t)dt\\
&=C_1\int_0^s t(-H'_j(t))dt+C_1sH_j(s)\\
&=C_1\int_0^sH_j(t)dt,
\end{align*}
\begin{align*}
(II)&=C_2(1-e^{-\theta s})\int_0^1 F_{\mu_j}(t)dt\\
&=C_2(1-e^{-\theta s})H_j(0),
\end{align*}
\begin{align*}
(III)&=C_2\int_0^1e^{-\theta t} F_{\mu_j}(t)dt\\
&=C_2\int_0^1e^{-\theta t}(-H'_j(t))dt\\
&=C_2H_j(0)-C_2\theta\int_0^1e^{-\theta t}H_j(t)dt
\end{align*}
and 
\[
-(III)=-C_2H_j(0)+C_2\theta\int_0^1e^{-\theta t}H_j(t)dt.
\]
Now,
\begin{align*}
(IV)&=C_2e^{\theta s}\int_s^1e^{-\theta t} F_{\mu_j}(t)dt\\
&=C_2e^{\theta s}\int_s^1e^{-\theta t}(-H'_j(t))dt\\
&=C_2e^{\theta s}\Big(e^{-\theta s}H_j(s)-\theta\int_s^1e^{-\theta t}H_j(t)dt\Big)\\
&=C_2H_j(s)-C_2\theta e^{\theta s}\int_s^1e^{-\theta t}H_j(t)dt
\end{align*}
and, finally,
\begin{align*}
(V)&=C_2e^{-\theta s}\int_0^se^{\theta t} F_{\mu_j}(t)dt\\
&=C_2e^{-\theta s}\int_0^se^{\theta t}(-H'_j(t))dt\\
&=C_2e^{-\theta s}\big(H_j(0)-e^{\theta s}H_j(s)\big)+C_2\theta e^{-\theta s}\int_0^se^{\theta t}H_j(t)dt\\
&=C_2e^{-\theta s}H_j(0)-C_2H_j(s)+C_2\theta e^{-\theta s}\int_0^se^{\theta t}H_j(t)dt.
\end{align*}
Therefore
\begin{align*}
\Gamma_\Delta F_{\mu_j}(s) &=(I)+(II)-(III)+(IV)+(V)\\&=C_1\int_0^sH_j(t)dt+C_2(1-e^{-\theta s})H_j(0)\\
&\hspace{0.5cm}-C_2H_j(0)+C_2\theta\int_0^1e^{-\theta t}H_j(t)dt+C_2H_j(s)-C_2\theta e^{\theta s}\int_s^1e^{-\theta t}H_j(t)dt\\
&\hspace{0.5cm}+C_2e^{-\theta s}H_j(0)-C_2H_j(s)+C_2\theta e^{-\theta s}\int_0^se^{\theta t}H_j(t)dt\\
&=C_1\int_0^sH_j(t)dt+C_2\theta\int_0^1e^{-\theta t}H_j(t)dt-C_2\theta e^{\theta s}\int_s^1e^{-\theta t}H_j(t)dt\\
&\hspace{0.5cm}+C_2\theta e^{-\theta s}\int_0^se^{\theta t}H_j(t)dt.
\end{align*}
Setting
\[
U(s)=e^{-\theta s}\int_0^se^{\theta t}H_j(t)dt-e^{\theta s}\int_s^1e^{-\theta t}H_j(t)dt,\quad s\in [0,1],
\]
we finally obtain that for any $s\in[0,1]$,
\begin{equation}\label{lambdaU}
\Gamma_\Delta F_{\mu_j}(s)= C_1\int_0^sH_j(t)dt+C_2\theta\int_0^1e^{-\theta t}H_j(t)dt+C_2\theta U(s).
\end{equation}
We then obtain that 
\begin{equation}\label{vpHawkes}
\Gamma_\Delta F_{\mu_j}(s)=\lambda_j F_{\mu_j}(s),\quad s\in[0,1]
\end{equation}
if and only if
\[
\lambda_j F_{\mu_j}(s)=C_1\int_0^sH_j(t)dt+C_2\theta\int_0^1e^{-\theta t}H_j(t)dt+C_2\theta U(s), \quad s\in[0,1].
\]
In the sequel, we use the following lemma.
\begin{lemma}
Setting
\begin{equation}\label{def-W}
Z(s)=\int_0^1 e^{-\theta|t-s|}H_j(t)dt,
\end{equation}
we have for any $s\in[0,1]$,
\[
U(s)=U(0)+2\int_0^sH_j(t)dt-\theta\int_0^sZ(t)dt.
\]
\end{lemma}
\begin{proof}[Proof of the lemma]
We have for any $s\in[0,1]$,
\begin{align*}
U'(s)&=-\theta e^{-\theta s}\int_0^se^{\theta t}H_j(t)dt+H_j(s)-\theta e^{\theta s}\int_s^1e^{-\theta t}H_j(t)dt+H_j(s)\\
&=2H_j(s)-\theta Z(s).
\end{align*}
This provides the result of the lemma.
\end{proof}

Now, we assume that \eqref{vpHawkes} is true. We have
\begin{align*}
\lambda_jF_{\mu_j}'(s)&=C_1H_j(s)+C_2\theta U'(s)\\
&=C_1H_j(s)+C_2\theta\Big(2H_j(s)-\theta Z(s)\Big)\\
&=\frac{C_1\theta^2}{\beta^2}H_j(s)-C_2\theta^2 Z(s),
\end{align*}
by using \eqref{C1C2} and the lemma.
This gives
\begin{equation}\label{F-W}
-\lambda_j H''_j(s)=\frac{C_1\theta^2}{\beta^2}H_j(s)-C_2\theta^2 Z(s), \quad s\in [0,1],
\end{equation}
meaning that $(\lambda_j,H_j)$ is solution of Equation \eqref{SL-Hawkes} of the main file.

\bigskip

Conversely, assume that $(\lambda_j,H_j)$ is solution of Equation \eqref{SL-Hawkes} of the main file. We prove  that \eqref{vpHawkes} is satisfied. Indeed, using \eqref{lambdaU}, the result of the lemma and
\[
U(0)=-\int_0^1e^{-\theta t}H_j(t)dt,
\]
 we have for any $s\in[0,1]$,
\begin{align*}
\Gamma_\Delta F_{\mu_j}(s)&= C_1\int_0^sH_j(t)dt+C_2\theta\int_0^1e^{-\theta t}H_j(t)dt+C_2\theta U(s)\\
&=C_1\int_0^sH_j(t)dt+C_2\theta\int_0^1e^{-\theta t}H_j(t)dt+C_2\theta\Big(U(0)+2\int_0^sH_j(t)dt-\theta\int_0^sZ(t)dt\Big)\\
&=\frac{C_1\theta^2}{\beta^2}\int_0^sH_j(t)dt-C_2\theta^2\int_0^sZ(t)dt.
\end{align*}
by using \eqref{C1C2}. Since $(\lambda_j,H_j)$ is solution of Equation \eqref{SL-Hawkes} of the main file, we finally obtain for any $s\in[0,1]$,
\begin{align*}
\Gamma_\Delta F_{\mu_j}(s)&= -\lambda_j\int_0^s H''_j(t)dt\\
&=\lambda_j\big(H'_j(0)-H'_j(s)\big)\\
&=\lambda_j F_{\mu_j}(s)
\end{align*}
and the result is proved.

\bigskip

Now, we prove the result of Remark~\ref{eigenvalue=0Hawkes} of the main file. Using Equation~\eqref{F-W}, we obtain
\[
-\lambda_j H^{(4)}_j(s)=\frac{C_1\theta^2}{\beta^2}H''_j(s)-C_2\theta^2 Z''(s), \quad s\in [0,1].
\]
By definition of $Z$ given in \eqref{def-W}, we have for any $s\in [0,1],$
\[
Z'(s)=-\theta e^{-\theta s}\int_0^se^{\theta u}H_j(u)du+\theta e^{\theta s}\int_s^1 e^{-\theta u}H_j(u)du
\]
and
\begin{align*}
Z''(s)&=\theta^2 e^{-\theta s}\int_0^se^{\theta u}H_j(u)du+\theta^2 e^{\theta s}\int_s^1 e^{-\theta u}H_j(u)du-2\theta H_j(s)\\
&=\theta^2Z(s)-2\theta H_j(s).
\end{align*}
Now, we obtain, still using Equations~\eqref{F-W} and \eqref{C1C2},
\begin{align*}
-\lambda_j H^{(4)}_j(s)&=\frac{C_1\theta^2}{\beta^2}H''_j(s)-C_2\theta^4 Z(s)+2C_2\theta^3H_j(s)\\
&=\frac{C_1\theta^2}{\beta^2}H''_j(s)-\lambda_j\theta^2H''_j(s)+\Big(2C_2\theta^3-\frac{C_1\theta^4}{\beta^2}\Big)H_j(s)\\
&=\frac{C_1\theta^2}{\beta^2}H''_j(s)-\lambda_j\theta^2H''_j(s)-C_1\theta^2H_j(s).
\end{align*}
Finally, we obtain for any $s\in [0,1],$
\begin{equation}\label{eqordre4}
-\lambda_j H^{(4)}_j(s)+\lambda_j\theta^2H''_j(s)=\frac{C_1\theta^2}{\beta^2}\Big(H''_j(s)-\beta^2H_j(s)\Big).
\end{equation}
In particular, if $\lambda_j=0$, then
\[
H_j(s)=c_1e^{\beta s}+c_2e^{-\beta s},\quad s\in [0,1],
\]
for $c_1$ and $c_2$ two constants.
The boundary condition $H_j'(0)=0$ gives $c_1=c_2$ and $H_j(1)=0$ gives $c_1=c_2=0$. Then $\lambda_j=0$ entails $H_j=0$ and then $F_{\mu_j}=0$. We have assumed that $F_{\mu_j}$ is continuous, which allows to use previous computations.
\subsection{Proof of Theorem~\ref{Hawkes-ed-sol}}
We set 
\[w_1=\frac{\beta w_0}{\beta-\alpha},\quad c=\frac{w_0\alpha\beta(2\beta-\alpha)}{2(\beta-\alpha)^2},\quad \theta=\beta-\alpha\]
and
\[
I(y)(t)=\frac{c}{w_1}\int_0^1e^{-\theta |t-s|}y(s)ds,\quad t\in[0,1].
\]
In the sequel, we assume that Assumption \eqref{cond-rho} of the main file is satisfied, which is equivalent to
\begin{equation}\label{cond-constantes1}
 \frac{2c}{w_1}\Big(2+\frac{3-2e^{-\theta/2}-e^{-\theta}}{\theta}\Big)<1.
\end{equation}
Now \eqref{SL-Hawkes} of the main file writes
\[
\left\{
\begin{array}{ll}
\displaystyle -\lambda y''(t)=w_1y(t)+w_1I(y(t)),   &t\in(0,1),\\
y(1)  =0,\quad y'(0)=0.&
\end{array}
\right.
\]
Let $H_j$ be a solution of this ODE as expressed in Theorem~\ref{Hawkes-ed}. We have
\[
\left\{
\begin{array}{ll}
\displaystyle H_j''(t)+\rho_j^2H_j(t)=-\rho_j^2I(H_j)(t),   &t\in(0,1),\\
H_j(1)  =0,\quad H_j'(0)=0&
\end{array}
\right.
\]
with
\begin{equation}\label{rhodef}
\rho_j:=\sqrt{\frac{w_1}{\lambda_j}}\underset{j\to+\infty}{\longrightarrow}+\infty.
\end{equation}
We write $H_j$ and $H_j'$ as 
\[
\left\{
\begin{array}{rc}
H_j(t)=&a(t)\cos (\rho_j t)+b(t)\sin (\rho_j t),\\
H'_j(t)=&-\rho_ja(t)\sin (\rho_j t)+\rho_jb(t)\cos (\rho_j t),
\end{array}
\right.
\]
with $a$ and $b$ $C^2$-functions that must satisfy
\[
\left\{
\begin{array}{rl}
a'(t)\cos (\rho_j t)+b'(t)\sin (\rho_j t)&=0\\
-a'(t)\rho_j\sin (\rho_j t)+b'(t)\rho_j\cos (\rho_j t)&=-\rho_j^2I(H_j)(t).
\end{array}
\right.
\]
This means that
\[
\left(\begin{array}{cc}\cos (\rho_j t) &\sin (\rho_j t) \\-\sin (\rho_j t) & \cos (\rho_j t)\end{array}\right)\left(\begin{array}{c}a'(t) \\ b'(t)\end{array}\right)=\left(\begin{array}{c}0 \\-\rho_jI(H_j)(t)\end{array}\right),
\]
or equivalently
\[
\left(\begin{array}{c}a'(t) \\ b'(t)\end{array}\right)=\left(\begin{array}{cc}\cos (\rho_j t) &-\sin (\rho_j t) \\\sin (\rho_j t) & \cos (\rho_j t)\end{array}\right)\left(\begin{array}{c}0 \\-\rho_jI(H_j)(t)\end{array}\right).
\]
Therefore
\[
\left\{
\begin{array}{rl}
\displaystyle a(t)&=a+\rho_j\int_0^t\sin (\rho_j s)I(H_j)(s)ds\\
\displaystyle b(t)&=b-\rho_j\int_0^t\cos (\rho_j s)I(H_j)(s)ds,
\end{array}
\right.
\]
with $a$ and $b$ two constants. 
Finally, for any $t\in [0,1],$
\[
H_j(t)=a\cos (\rho_j t)+b\sin (\rho_j t)+\rho_j\cos (\rho_j t)\int_0^t\sin (\rho_j s)I(H_j)(s)ds-\rho_j\sin (\rho_j t)\int_0^t\cos (\rho_j s)I(H_j)(s)ds,
\]
and
\[
H'_j(t)=-a\rho_j\sin (\rho_j t)+b\rho_j\cos (\rho_j t)-\rho_j^2\sin (\rho_j t)\int_0^t\sin (\rho_j s)I(H_j)(s)ds-\rho_j^2\cos (\rho_j t)\int_0^t\cos (\rho_j s)I(H_j)(s)ds,
\]
with $a$ and $b$ two constants. The boundary conditions $H'_j(0)=0$ implies that $b=0$. We then obtain that $H_j$ must write
\begin{equation}\label{Fjmust}
H_j(t)=a\cos (\rho_j t)+\rho_j\cos (\rho_j t)\int_0^t\sin (\rho_j s)I(H_j)(s)ds-\rho_j\sin (\rho_j t)\int_0^t\cos (\rho_j s)I(H_j)(s)ds,
\end{equation}
and
\[
H'_j(t)=-a\rho_j\sin (\rho_j t)-\rho_j^2\sin (\rho_j t)\int_0^t\sin (\rho_j s)I(H_j)(s)ds-\rho_j^2\cos (\rho_j t)\int_0^t\cos (\rho_j s)I(H_j)(s)ds.
\]
In particular, we have
\begin{align*}
H_j(t)&=a\cos (\rho_j t)-\rho_j\int_0^t\sin (\rho_j (t-s))I(H_j)(s)ds\\
&=a\cos(\rho_j t)-\int_0^{\rho_j t}\sin(u)I\big(H_j\big)\Big(t-\frac{u}{\rho_j}\Big)du.
\end{align*}
We denote $J_j:{\mathbb L}_\infty[0,1]\mapsto {\mathbb L}_\infty[0,1]$ the linear operator defined for $f\in {\mathbb L}_\infty[0,1]$, for $t\in [0,1]$,  by
\begin{align*}
J_j(f)(t)&:=\rho_j\cos (\rho_j t)\int_0^t\sin (\rho_j s)I(f)(s)ds-\rho_j\sin (\rho_j t)\int_0^t\cos (\rho_j s)I(f)(s)ds\\
&=-\rho_j\int_0^t\sin (\rho_j (t-s))I(f)(s)ds\\
&=-\int_0^{\rho_j t}\sin(u)I\big(f\big)\Big(t-\frac{u}{\rho_j}\Big)du.
\end {align*}
Finally, if $H_j$ satisfies Equation~\eqref{SL-Hawkes} of the main file, $H_j$ satisfies
\begin{equation}\label{FJ}
H_j(t)=a\cos (\rho_j t)+J_j(H_j)(t),\quad t\in [0,1].
\end{equation}
Conversely, if $H_j$ satisfies Equation~\eqref{FJ}, for $t\in [0,1]$, 
\begin{align*}
H''_j(t)&=-a\rho_j^2\cos (\rho_j t)+J''_j(H_j)(t)\\
&=-a\rho_j^2\cos (\rho_j t)-\rho_j^3\cos (\rho_j t)\int_0^t\sin (\rho_j s)I(H_j)(s)ds-\rho_j^2\sin^2(\rho_j t)I(H_j)(t)\\
&\hspace{1cm}+\rho_j^3\sin (\rho_j t)\int_0^t\cos (\rho_j s)I(H_j)(s)ds-\rho_j^2\cos^2 (\rho_j t)I(H_j)(t)\\
&=-a\rho_j^2\cos (\rho_j t)-\rho_j^2J_j(H_j)(t)-\rho_j^2I(H_j)(t)\\
&=-\rho_j^2H_j(t)-\rho_j^2I(H_j)(t).
\end {align*}
It means that if $H_j$ satisfies Equation~\eqref{FJ}, then it satisfies \eqref{SL-Hawkes} of the main file. We have proved:
\begin{equation}\label{equiv-equa}
H_j \mbox{ satisfies Equation \eqref{SL-Hawkes} of the main file} \iff H_j \mbox{ satisfies Equation  \eqref{FJ}} .
\end{equation}

Now, let us deal with $F_{\mu_j}(t)=-H_j'(t)$, for $t\in [0,1]$. We have:
\begin{align*}
F_{\mu_j}(t)&=-H_j'(t)\\
&=a\rho_j\sin (\rho_j t)+\rho_j^2\sin (\rho_j t)\int_0^t\sin (\rho_j s)I(H_j)(s)ds+\rho_j^2\cos (\rho_j t)\int_0^t\cos (\rho_j s)I(H_j)(s)ds\\
&=a\rho_j\sin (\rho_j t)+\rho_j^2\int_0^t\cos(\rho_j (t-s))I(H_j)(s)ds\\
&=a\rho_j\sin (\rho_j t)+\rho_j\int_0^{\rho_j t}\cos(u)I\big(H_j\big)\Big(t-\frac{u}{\rho_j}\Big)du\\
&=a\rho_j\sin (\rho_j t)+\frac{c\rho_j}{w_1}\int_0^{\rho_j t}\cos(u)\Big(\int_0^1e^{-\theta\left|t-\rho_j^{-1}u-s\right|}H_j(s)ds\Big)du.\end{align*}
It gives
\begin{equation}\label{KJ}
F_{\mu_j}(t)=a\rho_j\sin (\rho_j t)+\rho_jK_j(H_j)(t),\quad t\in [0,1],
\end{equation}
with for any bounded function $f$,
\[
K_j(f)(t)=\frac{c}{w_1}\int_0^{\rho_j t}\cos(u)\Big(\int_0^1e^{-\theta\left|t-\rho_j^{-1}u-s\right|}f(s)ds\Big)du.
\]
\begin{lemma}\label{Opsup}
We set
\[
\Vvert J_j \Vvert_{\infty}:=\sup_{f:\ \|f\|_{\infty}=1} \|J_j(f)\|_{\infty}.
\]
If \eqref{cond-constantes1} is satisfied, we have 
\[
\Vvert J_j\Vvert_{\infty}<\frac{1}{2}.
\]
Furthermore, there exist $c_1$, $c_2$ and $c_3$ three positive constants only depending on $c$, $\theta$ and $w_1$ such that
\[
\|J_j(\cos(\rho_j\cdot))\|_{\infty}\leq c_1\rho_j^{-1},\quad
\|K_j(\cos(\rho_j\cdot))\|_{\infty}\leq c_2\rho_j^{-1}\quad
\mbox{and}\quad \Vvert K_j\Vvert_{\infty}\leq c_3.
\]
\end{lemma}
\begin{proof}[Proof of Lemma~\ref{Opsup}]
 We have
\begin{align*}
-J_j(f)(t)&=\int_0^{\rho_j t}\sin(u)I\big(f\big)\Big(t-\frac{u}{\rho_j}\Big)du\\
&=\frac{c}{w_1}\int_0^{\rho_j t}\sin(u)\Big(\int_0^1e^{-\theta\left|t-\rho_j^{-1}u-s\right|}f(s)ds\Big)du\\
&=:{\mathcal I}m\Bigg(J^{(1)}(f)(t)+J^{(2)}(f)(t)+J^{(3)}(f)(t)\Bigg),
\end{align*}
with
\begin{align*}
J^{(1)}(f)(t)&=\frac{c}{w_1}\int_0^{\rho_j t}e^{iu}\Big(\int_0^{t-\rho_j^{-1}u}e^{-\theta\left(t-\rho_j^{-1}u-s\right)}f(s)ds\Big)du\\
&=\frac{c}{w_1}\int_0^te^{-\theta(t-s)}f(s)\Big(\int_0^{\rho_j (t-s)}e^{(i+\theta\rho_j^{-1})u}du\Big)ds\\
&=\frac{c}{w_1(i+\theta\rho_j^{-1})}\int_0^te^{-\theta(t-s)}f(s)\Big(e^{i\rho_j (t-s)+\theta(t-s)}-1\Big)ds\\
&=\frac{c(-i+\theta\rho_j^{-1})}{w_1\sqrt{1+\theta^2\rho_j^{-2}}}\int_0^tf(s)\Big(e^{i\rho_j (t-s)}-e^{-\theta(t-s)}\Big)ds,
\end{align*}
\begin{align*}
J^{(2)}(f)(t)&=\frac{c}{w_1}\int_0^{\rho_j t}e^{iu}\Big(\int_{t-\rho_j^{-1}u}^te^{\theta\left(t-\rho_j^{-1}u-s\right)}f(s)ds\Big)du\\
&=\frac{c}{w_1}\int_0^te^{\theta(t-s)}f(s)\Big(\int_{\rho_j (t-s)}^{\rho_j t}e^{(i-\theta\rho_j^{-1})u}du\Big)ds\\
&=\frac{c}{w_1(i-\theta\rho_j^{-1})}\int_0^te^{\theta(t-s)}f(s)\Big(e^{i\rho_jt-\theta t}-e^{i\rho_j(t-s)-\theta (t-s)}\Big)ds\\
&=\frac{c(-i-\theta\rho_j^{-1})}{w_1\sqrt{1+\theta^2\rho_j^{-2}}}\int_0^tf(s)\Big(e^{i\rho_jt-\theta s}-e^{i\rho_j(t-s)}\Big)ds
\end{align*}
and
\begin{align*}
J^{(3)}(f)(t)&=\frac{c}{w_1}\int_0^{\rho_j t}e^{iu}\Big(\int_t^1e^{\theta\left(t-\rho_j^{-1}u-s\right)}f(s)ds\Big)du\\
&=\frac{c}{w_1}\int_t^1e^{\theta(t-s)}f(s)\Big(\int_{0}^{\rho_j t}e^{(i-\theta\rho_j^{-1})u}du\Big)ds\\
&=\frac{c}{w_1(i-\theta\rho_j^{-1})}\int_t^1e^{\theta(t-s)}f(s)\Big(e^{i\rho_jt-\theta t}-1\Big)ds\\
&=\frac{c(-i-\theta\rho_j^{-1})}{w_1\sqrt{1+\theta^2\rho_j^{-2}}}\int_t^1f(s)\Big(e^{i\rho_jt-\theta s}-e^{\theta(t-s)}\Big)ds.
\end{align*}
This yields
\begin{align}\label{majoJj}
-J_j(f)(t)&=\frac{2c\theta\rho_j^{-1}}{w_1\sqrt{1+\theta^2\rho_j^{-2}}}{\mathcal I}m\Bigg(\int_0^tf(s)e^{i\rho_j (t-s)}ds\Bigg)+\frac{c}{w_1\sqrt{1+\theta^2\rho_j^{-2}}}\int_0^1f(s)e^{-\theta|t-s|}ds\\
&\hspace{1cm}+{\mathcal I}m\Bigg(\frac{c(-i-\theta\rho_j^{-1})}{w_1\sqrt{1+\theta^2\rho_j^{-2}}}\int_0^1f(s)e^{i\rho_jt-\theta s}ds\Bigg).\nonumber
\end{align}
We obtain for any $t\in [0,1]$,
\begin{align*}
|J_j(f)(t)|&\leq\Big(\frac{2c}{w_1}+\frac{c}{w_1}\int_0^1e^{-\theta|t-s|}ds+\frac{c}{w_1}\int_0^1e^{-\theta s}ds\Big)\|f\|_{\infty}\\
&\leq\Big(\frac{2c}{w_1}+\frac{ce^{-\theta t}}{w_1}\int_0^te^{\theta s}ds+\frac{ce^{\theta t}}{w_1}\int_t^1e^{-\theta s}ds+\frac{c}{w_1}\int_0^1e^{-\theta s}ds\Big)\|f\|_{\infty}\\
&\leq\frac{c}{w_1}\Big(2+\frac{3-e^{-\theta t}-e^{\theta(t-1)}-e^{-\theta}}{\theta}\Big)\|f\|_{\infty}\\
&\leq\frac{c}{w_1}\Big(2+\frac{3-2e^{-\theta/2}-e^{-\theta}}{\theta}\Big)\|f\|_{\infty}
\end{align*}
and
\[
\Vvert J_j\Vvert_{\infty}\leq \frac{c}{w_1}\Big(2+\frac{3-2e^{-\theta/2}-e^{-\theta}}{\theta}\Big).
\]
This yields the first result.

For the second result, observe that with $f=\cos(\rho_j\cdot)$, we obtain from \eqref{majoJj},
\begin{align*}
|J_j(\cos(\rho_j\cdot))(t)|&\leq\frac{2c\theta}{w_1\rho_j}+\frac{c}{w_1}\Big|\int_0^1e^{i\rho_j s}e^{-\theta |t-s|}ds\Big|+\frac{c}{w_1}\Big|\int_0^1e^{i\rho_j s}e^{-\theta s}ds\Big|.
\end{align*}
But
\begin{align*}
\int_0^1e^{i\rho_j s}e^{-\theta |t-s|}ds&=e^{-\theta t}\int_0^te^{(i\rho_j+\theta) s}ds+e^{\theta t}\int_t^1e^{(i\rho_j-\theta) s}ds\\
&=e^{-\theta t}\Big(\frac{e^{(i\rho_j+\theta)t}-1}{i\rho_j+\theta}\Big)+e^{\theta t}\Big(\frac{e^{i\rho_j-\theta}-e^{(i\rho_j-\theta)t}}{i\rho_j-\theta}\Big)
\end{align*}
and there exists $c_1$ only depending on $c$, $\theta$ and $w_1$ such that
\[
\|J_j(\cos(\rho_j\cdot))\|_{\infty}\leq c_1\rho_j^{-1}.
\]
For the third result, observe that
\begin{align*}
K_j(f)(t)&={\mathcal R}e\Bigg(\frac{c}{w_1}\int_0^{\rho_j t}e^{iu}\Big(\int_0^1e^{-\theta\left|t-\rho_j^{-1}u-s\right|}f(s)ds\Big)du\Bigg)\\
&={\mathcal R}e\Bigg(J^{(1)}(f)(t)+J^{(2)}(f)(t)+J^{(3)}(f)(t)\Bigg)\\
&=\frac{2c\theta\rho_j^{-1}}{w_1\sqrt{1+\theta^2\rho_j^{-2}}}{\mathcal R}e\Bigg(\int_0^tf(s)e^{i\rho_j (t-s)}ds\Bigg)+\frac{c\theta\rho_j^{-1}}{w_1\sqrt{1+\theta^2\rho_j^{-2}}}\int_0^1f(s)e^{-\theta|t-s|}ds\\
&\hspace{1cm}+{\mathcal R}e\Bigg(\frac{c(-i-\theta\rho_j^{-1})}{w_1\sqrt{1+\theta^2\rho_j^{-2}}}\int_0^1f(s)e^{i\rho_jt-\theta s}ds\Bigg).
\end{align*}
Therefore, as previously, there exists $c_2$ only depending on $c$, $\theta$ and $w_1$ such that
\[
\|K_j(\cos(\rho_j\cdot))\|_{\infty}\leq c_2\rho_j^{-1}.
\]
and  there exists $c_3$ only depending on $c$, $\theta$ and $w_1$ such that 
\[
\Vvert K_j\Vvert_{\infty}\leq c_3.
\]
\end{proof}
Now, we end the proof of the theorem and use Equivalence~\eqref{equiv-equa}. Since $J_j$ is a linear operator, for any integer $N$, we have, for $t\in[0,1]$,
\begin{align*}
H_j(t)&=a\cos (\rho_j t)+J_j(H_j)(t)\\
&=a\cos (\rho_j t)+aJ_j(\cos(\rho_j\cdot))(t)+(J_j\circ J_j)(H_j)(t)\\
&=a\sum_{k=0}^{N-1}J_j^{k\circ}(\cos(\rho_j\cdot))(t)+(J_j^{N\circ})(H_j)(t),
\end{align*}
with for any bounded function $f$,
\[
J_j^{k\circ}(f)=(\underbrace{J_j\circ\cdots \circ J_j}_{k \textrm{ times}})(f).
\]
Since $H_j$ is bounded and $\Vvert J_j\Vvert_{\infty}<1/2$ (see Lemma~\ref{Opsup}), we have
\[
\limsup_{N\to +\infty}\|(J_j^{N\circ})(H_j)\|_{\infty}\leq \limsup_{N\to +\infty}\Vvert J_j\Vvert_{\infty}^N\|H_j\|_{\infty}=0.
\]
Therefore
\begin{equation}\label{Fjexpan}
H_j(t)=a\sum_{k=0}^{+\infty}J_j^{k\circ}(\cos(\rho_j\cdot))(t),\quad t\in [0,1].
\end{equation}
This implies that, up to renormalization, we have at most one solution to the problem.
Now, we obtain a solution of the system \eqref{SL-Hawkes} of the main file if and only if the boundary condition $H_j(1)=0$ is satisfied. Using \eqref{Fjexpan}, this is equivalent to
\[
\cos(\rho_j)=-\sum_{k=1}^{+\infty}J_j^{k\circ}(\cos(\rho_j\cdot))(1).
\]
We can write $\rho_j=(k_j-1/2)\pi+u_j\pi$, with $k_j\in{\mathbb N}^*$ and $-1/2<u_j\leq1/2$ ($k_j$ and $u_j$ are uniquely defined). Then, $H_j(1)=0$ writes
\begin{equation}\label{uj}
(-1)^{k_j}\sin(u_j\pi)+\sum_{k=1}^{+\infty}J_j^{k\circ}(\cos(((k_j-1/2)\pi+u_j\pi)\cdot))(1)=0.
\end{equation}
Since $\Vvert  J_j\Vvert_{\infty}<1/2$ (see Lemma~\ref{Opsup}), for $k_j$ fixed, the left hand side is a continuous fonction of $u_j$ taking positive and negative values when $u_j$ describes the interval $(-1/2;1/2]$. Therefore, for any positive integer $k_j$, there exists $u_j\in (-1/2;1/2]$ such that Equation~\eqref{uj} is satisfied. Now, let us give a control of $u_j$. We have:
\begin{align*}
2|u_j|\leq|\sin(u_j\pi)|&\leq \sum_{k=1}^{+\infty}\Vvert J_j\Vvert_{\infty}^{k-1}\times\|J_j(\cos(\rho_j\cdot))\|_{\infty}\leq \sum_{k=1}^{+\infty}\Big(\frac{1}{2}\Big)^{k-1}\times\frac{c_1}{\rho_j}=\frac{2c_1}{\rho_j}.
\end{align*}
leading to
\[
|u_j|\leq \frac{c_1}{\rho_j}.
\]
These computations hold for any integer $k_j$ and in particular, we can consider $k_j=j$. In this case, we obtain
\[
\rho_j=(j-1/2)\pi+u_j\pi
\]
and, using \eqref{rhodef},
\[
\lambda_j=\frac{w_1}{\rho_j^2}=\frac{w_1}{\big((j-1/2)\pi+u_j\pi)\big)^2}=\frac{w_1}{\big(j\pi-\pi/2\big)^2}\Big(1+\frac{u_j}{(j-1/2)}\Big)^{-2}.
\]
We have finally proved that for any $j\geq 1$, System \eqref{SL-Hawkes} of the main file has a solution $(\lambda_j,H_j)$ and there exists a constant $C_1$ depending only on $c$, $\theta$ and $w_1$ such that for any $j\geq 1$, 
\[
\Big|\lambda_j-\frac{w_1}{\big(j\pi-\pi/2\big)^2}\Big|\leq C_1j^{-4}.
\]
Now, it remains to determine $F_{\mu_j}$. Using \eqref{KJ}, we have for any $t\in [0,1]$,
\begin{align*}
F_{\mu_j}(t)&=a\rho_j\sin (\rho_j t)+\rho_jK_j(H_j)(t)\\
&=a\rho_j\sin (\rho_j t)+a\rho_jK_j(\cos(\rho_j\cdot))(t)+a\rho_j\sum_{k=1}^{+\infty}K_j\big(J_j^{k\circ}(\cos(\rho_j\cdot))\big)(t).
\end{align*}
and we have to determine the value of $a$ by using the property $\|F_{\mu_j}\|_2=1$. Setting 
\[
v_j(t):=\rho_jK_j(\cos(\rho_j\cdot))(t)+\rho_j\sum_{k=1}^{+\infty}K_j\big(J_j^{k\circ}(\cos(\rho_j\cdot))\big)(t),\quad t\in[0,1],
\]
so that
\begin{equation}\label{Fjvj}
F_{\mu_j}(t)=a\rho_j\sin (\rho_j t)+av_j(t),\quad t\in[0,1],
\end{equation}
we have, using Lemma~\ref{Opsup},
\[
\|v_j\|_{\infty}\leq c_2+c_3\rho_j\sum_{k=1}^{+\infty}2^{1-k}\|J_j(\cos(\rho_j\cdot))\|_{\infty}\leq c_2+2c_1c_3
\]
and, since $\rho_j^{-1}= O(j^{-1})$,
\begin{align*}
1=\int F^2_{\mu_j}(t)dt&=a^2\rho_j^2\int_0^1\Big(\sin(\rho_j t)+\frac{v_j(t)}{\rho_j}\Big)^2dt=a^2\rho_j^2\times\Big(\frac{1}{2}+O(j^{-1})\Big).
\end{align*}
This implies $a\rho_j=\sqrt{2}+O(j^{-1})$ and $a=\frac{\sqrt{2}}{\rho_j}+O(j^{-2})$.
Using again Equation~\eqref{Fjvj} yields the result stated in Equation~\eqref{etaj} of the main file. Equation~\eqref{Fjexpan} shows that the eigenvalues $(\lambda_j)_{j\geq 1} $ are simple.
\section{Proofs of Section~\ref{sec:RMPC}}\label{sec:proof:op}
\subsection{Proof of Lemma~\ref{eq:construct_GDelta}}
We first prove some intermediate results. First remark that given $i=1,\hdots, n$, if $T\in N_i$, there exists $k\in\{0,\hdots,|\mathcal T|\}$ such that $T=T_k$, then for all $\ell\in \{1,\hdots,|\mathcal T|\}$, 
 \[
\int_{T_{\ell-1}}^{T_{\ell}}\mathbf 1_{T\leq t}dt=\int_{T_{\ell-1}}^{T_{\ell}}\mathbf 1_{T_{k}\leq t}dt=(T_{\ell}-T_{\ell-1})\mathbf 1_{T_k\leq T_{\ell-1}}=(T_{\ell}-T_{\ell-1})\mathbf 1_{T\leq T_{\ell-1}}. \]
This implies
\begin{eqnarray*}
\int_{T_{\ell-1}}^{T_\ell}F_{\Pi_i}(t)dt& =& \int_{T_{\ell-1}}^{T_\ell}\sum_{T\in N_i}\mathbf 1_{T\leq t}dt=\sum_{T\in N_i}\int_{T_{\ell-1}}^{T_\ell}\mathbf 1_{T\leq t}dt\\
&=&\sum_{T\in N_i}(T_\ell-T_{\ell-1})\mathbf 1_{T\leq T_{\ell-1}}=(T_\ell-T_{\ell-1})F_{\Pi_i}(T_{\ell-1}),
\end{eqnarray*}
and consequently
\begin{eqnarray*}
\int_{T_{\ell-1}}^{T_\ell}F_{\widehat W}(t)dt&=&(T_\ell-T_{\ell-1})\frac1n\sum_{i=1}^nF_{\Pi_i}(T_{\ell-1})  = (T_\ell-T_{\ell-1})F_{\widehat W}(T_{\ell-1}).  
\end{eqnarray*}
and also
\begin{equation}\label{intF}
\int_{T_{\ell-1}}^{T_\ell}F_{\widehat\Delta_i}(t)dt = (T_\ell-T_{\ell-1})F_{\widehat\Delta_i}(T_{\ell-1}).  
\end{equation}

We turn now to the proof of the first point. Let $\ell\in \{1,\hdots,|\mathcal T|\}$, 
 \begin{eqnarray*}
 \widehat\Gamma_{\widehat\Delta} e_\ell(s)& = &\frac1{\sqrt{T_{\ell}-T_{\ell-1}}}\int_{T_{\ell-1}}^{T_\ell}\widehat K_{\widehat\Delta}(s,t)dt\\
 &=&\frac1{\sqrt{T_{\ell}-T_{\ell-1}}}\int_{T_{\ell-1}}^{T_\ell}\left(\frac1n\sum_{i=1}^n\sum_{T,T'\in N_i}\mathbf 1_{T\leq s}\mathbf 1_{T'\leq t}-F_{\widehat W}(s)F_{\widehat W}(t)\right)dt\\
 &=&\frac1{\sqrt{T_{\ell}-T_{\ell-1}}}\frac1n\sum_{i=1}^n\left(\sum_{T,T'\in N_i}\mathbf 1_{T\leq s}\int_{T_{\ell-1}}^{T_\ell}\mathbf 1_{T'\leq t}dt-F_{\widehat W}(s)\int_{T_{\ell-1}}^{T_\ell}F_{\widehat W}(t)dt\right)\\
 &=&\frac1{\sqrt{T_{\ell}-T_{\ell-1}}}\frac1n\sum_{i=1}^n\left(F_{\Pi_i}(s)\sum_{T'\in N_i}(T_\ell-T_{\ell-1})\mathbf 1_{T'\leq T_{\ell-1}}-F_{\widehat W}(s)     F_{\widehat W}(T_{\ell-1})(T_\ell-T_{\ell-1})\right).\\
 &=&\frac{\sqrt{T_\ell-T_{\ell-1}}}n\sum_{i=1}^n\left(F_{\Pi_i}(s)F_{\Pi_i}(T_{\ell-1})-F_{\widehat W}(s)     F_{\widehat W}(T_{\ell-1})\right).
 \end{eqnarray*}
Following the same lines, for all $\ell'\in\{1,\hdots,|\mathcal T|\}$, 
\begin{eqnarray*}
\langle\widehat\Gamma_{\widehat\Delta} e_\ell,e_{\ell'}\rangle&=&\frac{1}{\sqrt{T_{{\ell'}}-T_{{\ell'}-1}}}\int_{T_{{\ell'}-1}}^{T_{{\ell'}}}\widehat\Gamma_{\widehat\Delta} e_\ell(s)ds\\
&=&\frac{\sqrt{T_\ell-T_{\ell-1}}}{n\sqrt{T_{{\ell'}}-T_{{\ell'}-1}}}\sum_{i=1}^n\int_{T_{{\ell'}-1}}^{T_{{\ell'}}}\left(F_{\Pi_i}(s)F_{\Pi_i}(T_{\ell-1})-F_{\widehat W}(s)     F_{\widehat W}(T_{\ell-1})\right)ds\\
&=&\frac{\sqrt{(T_\ell-T_{\ell-1})(T_{{\ell'}}-T_{{\ell'}-1})}}n\sum_{i=1}^n\left(F_{\Pi_i}(T_{\ell'-1})F_{\Pi_i}(T_{\ell-1})-F_{\widehat W}(T_{\ell'-1})     F_{\widehat W}(T_{\ell-1})\right).
\end{eqnarray*}
In particular, if for all $\ell$, there exists a unique $i$ such that $T_\ell\in N_i$, we have $F_{\widehat W}(T_{\ell-1})=(\ell-1)/n$ and
\[
\int_{T_{\ell-1}}^{T_\ell}F_{\widehat W}(t)dt=\frac{\ell-1}n(T_\ell-T_{\ell-1}),
\]
which concludes the proof of the first point.

Now, we start the proof of the second point by proving that 
\begin{equation}\label{eq:incluImST}
{\rm Im}(\widehat\Gamma_{\widehat\Delta})\subset{\rm span}\{e_1,\hdots,e_{|\mathcal T|} \}. 
\end{equation}
Since $\widehat\Gamma_{\widehat\Delta}$ is a finite rank operator, both spaces ${\rm Im}(\widehat\Gamma_{\widehat\Delta})$ and ${\rm span}\{e_1,\hdots,e_{|\mathcal T|}\}$ are finite-dimensional vector spaces and it is sufficient to prove 
\[
{\rm span}\{e_1,\hdots,e_{|\mathcal T|} \}^{\perp}\subset {\rm Im}(\widehat\Gamma_{\widehat\Delta})^{\perp}. 
\]
Moroever, since $\widehat\Gamma_{\widehat\Delta}$ is a self-adjoint operator, 
\[
 {\rm Im}(\widehat\Gamma_{\widehat\Delta})^{\perp}= {\rm Ker}(\widehat\Gamma_{\widehat\Delta}). 
 \]
 Then, to prove~\eqref{eq:incluImST}, it is sufficient to prove that 
 \[
{\rm span}\{e_1,\hdots,e_{|\mathcal T|} \}^{\perp}\subset  {\rm Ker}(\widehat\Gamma_{\widehat\Delta}). 
\]
Let $f\in{\rm span}\{e_1,\hdots,e_{|\mathcal T|} \}^{\perp}$, we have for any $\ell\in\{1,\ldots,|\mathcal T|\}$,
\[
\langle f,e_\ell\rangle =0,
\]
then, by definition of $e_\ell$, 
\[
\int_{T_{\ell-1}}^{T_\ell}f(t)dt=0. 
\]
This implies that, for all $i=1,\hdots,n$, for all $T\in N_i$, since, by construction of the grid $\mathcal T$,  there exists $\ell$ such that $T = T_{\ell-1}$, 
\[
\int_{T}^1 f(t) dt = \sum_{\ell'=\ell}^{|\mathcal T|}\int_{T_{\ell'}-1}^{T_{\ell'}}f(t)dt = 0. 
\]
Hence
\[
\langle F_{\widehat W}, f\rangle = \int_0^1 \widehat W([0,t]) f(t)dt =\frac1n\sum_{i=1}^n\sum_{T\in N_i}\int_0^1\mathbf 1_{\{T\leq t\}}f(t)dt = \frac1n\sum_{i=1}^n\sum_{T\in N_i}\int_T^1f(t)dt=0
\]
and 
\begin{eqnarray*}
\widehat\Gamma_{\widehat\Delta}f (t) &=&  \frac1n\sum_{i=1}^n\sum_{T,T'\in N_i}\int_0^1\mathbf 1_{\{T\leq s, T'\leq t\}}f(s) ds-\widehat W([0,t])\langle F_{\widehat W},f\rangle,\\
&=&\frac1n\sum_{i=1}^n\sum_{T,T'\in N_i}\mathbf 1_{\{T'\leq t\}}\int_T^1f(s) ds-\widehat W([0,t])\langle F_{\widehat W},f\rangle = 0. 
\end{eqnarray*}
This concludes the proof of \eqref{eq:incluImST}. 

Now we identify the eigenfunctions of $\widehat\Gamma_{\widehat\Delta}$ with the eigenvectors of the matrix $\widehat G_{\widehat\Delta}$. 
\begin{itemize}
\item Let $j\geq 1$. We consider $\widehat\lambda_j$ a non-zero eigenvalue of  $\widehat \Gamma_{\widehat\Delta}$. Let $\widehat\eta_j$ an associated eigenfunction such that $\|\widehat\eta_j\|=1$. We set
\[
\widehat v_j = \left(\langle \widehat\eta_j,e_\ell\rangle\right)_{\ell=1,\hdots,|\mathcal T|}. 
\]
We have 
\begin{eqnarray}
\widehat G_{\widehat\Delta}\widehat v_j &=&\left( \sum_{\ell'=1}^{|\mathcal T|}\langle\widehat\Gamma_{\widehat\Delta}e_{\ell},e_{\ell'}\rangle\langle \widehat\eta_j,e_{\ell'}\rangle\right)_{\ell=1,\hdots,|\mathcal T|}
\nonumber\\
&=&\left(\left\langle \widehat\Gamma_{\widehat\Delta}e_\ell,\sum_{\ell'=1}^{|\mathcal T|}\langle \widehat\eta_j,e_{\ell'}\rangle e_{\ell'}\right\rangle\right)_{\ell=1,\hdots,|\mathcal T|}.\label{eq:GDelta} \end{eqnarray}
Now remark that $\sum_{\ell'=1}^{|\mathcal T|}\langle \widehat\eta_j,e_{\ell'}\rangle e_{\ell'}$ is the orthogonal projection of the function $\widehat\eta_{j}$ on the space ${\rm span}\{e_1,\hdots,e_{|\mathcal T|}\}$. Then, since $\widehat\Gamma_{\widehat\Delta}\widehat\eta_j=\widehat\lambda_j\widehat\eta_j$  and $\widehat\lambda_j>0$, from \eqref{eq:incluImST}, we deduce  
\[
\widehat\eta_j\in{\rm Im}(\widehat\Gamma_{\widehat\Delta})\subset{\rm span}\{e_1,\hdots,e_{|\mathcal T|}\}, 
\]
meaning that 
\[
\sum_{\ell'=1}^{|\mathcal T|}\langle \widehat\eta_j,e_{\ell'}\rangle e_{\ell'} = \widehat\eta_j
\]
and $\|\widehat v_j \|=1$.
Then, since $\widehat\Gamma_{\widehat\Delta}$ is a self-adjoint operator, ~\eqref{eq:GDelta} can be written
\[
\widehat G_{\widehat\Delta}\widehat v_j = \left(\left\langle \widehat\Gamma_{\widehat\Delta}e_\ell,\widehat\eta_j\right\rangle\right)_{\ell=1,\hdots,|\mathcal T|}=\left(\left\langle e_\ell,\widehat\Gamma_{\widehat\Delta}\widehat\eta_j\right\rangle\right)_{\ell=1,\hdots,|\mathcal T|} = \widehat\lambda_j\left(\left\langle e_\ell,\widehat\eta_j\right\rangle\right)_{\ell=1,\hdots,|\mathcal T|}=\widehat\lambda_j \widehat v_j .
\]
This means that $\widehat\lambda_j\in {\rm Sp}(\widehat G_{\widehat\Delta})$ and, consequently 
\[
{\rm Sp}(\widehat\Gamma_{\widehat\Delta})\backslash\{0\}\subset {\rm Sp}(\widehat G_{\widehat\Delta}). 
\]
\item Conversely, let $\widehat v_j=(\widehat v_j^1,\hdots,\widehat v_j^{|\mathcal{T}|})^t$ a unit-norm eigenvector of $\widehat G_{\widehat\Delta}$ associated with the eigenvalue $\widehat\lambda_j$ and let
	$$
	\widehat\eta_j = \sum_{\ell=1}^{|\mathcal{T}|}\widehat v_j^\ell e_\ell. 
	$$
We have:	
\begin{eqnarray*}
\widehat\Gamma_{\widehat\Delta}\widehat\eta_j  &=& \sum_{\ell'=1}^{|\mathcal T|}\langle\widehat\Gamma_{\widehat\Delta}\widehat\eta_j , e_{\ell'}\rangle e_{\ell'} =\sum_{\ell,\ell'=1}^{|\mathcal T|}\widehat v_j^\ell \langle\widehat\Gamma_{\widehat\Delta}e_\ell,e_{\ell'}\rangle e_{\ell'} = \sum_{\ell'=1}^{|\mathcal T|}\left[\widehat G_{\widehat\Delta}\widehat v_j\right]_{\ell'}e_{\ell'} \\
&=&\widehat\lambda_j\sum_{\ell'=1}^{|\mathcal T|}\widehat v_j^{\ell'} e_{\ell'} = \widehat\lambda_j\widehat\eta_j . 
\end{eqnarray*}
	Then, $\widehat\eta_j$ is a unit-norm eigenfunction of the operator $\widehat\Gamma_{\widehat\Delta}$ associated with the eigenvalue $\widehat\lambda_j$
 and 
\[{\rm Sp}(\widehat G_{\widehat\Delta})\subset
{\rm Sp}(\widehat\Gamma_{\widehat\Delta}). 
\]
To end the proof of the third point, it remains to remark that, since the functions $e_\ell$'s are right continuous functions, then $\widehat\eta_j$ also is a right-continuous function. Moreover, following the same lines than the proof of Proposition~\ref{prop:muj}, $\widehat\eta_j$ is a function with bounded variations and $\widehat\eta_j(t)=0$ for $t<0$. Then, applying Proposition~4.4.3 of \cite{Cohn93}, we obtain the expected result: there exists a unique measure $\widehat\mu_j$ such that  $F_{\widehat\mu_j} = \widehat\eta_j. $
\end{itemize}
%


\subsection{Proof of Theorem~\ref{thm:CV_rates}}
We first prove the first point. From \cite{Bosq}, Lemma 4.2, p.~103, we can deduce that 
\begin{eqnarray}\label{eq:Bosqeiv}
\sup_{j\geq 1}|\widehat\lambda_j-\lambda_j|\leq \Vvert\Gamma_\Delta-\widehat\Gamma_{\widehat\Delta}\Vvert,
\end{eqnarray}
where, for an operator $S$ on $\mathbb L^2([0,1])$,  $$\Vvert S\Vvert=\sup_{x\in\mathbb L^2([0,1])}\frac{\|Sx\|}{\|x\|}$$ denotes the usual operator norm.
Let  $\widehat\Gamma_\Delta = \frac1n\sum_{i=1}^nF_{\Delta_i}\otimes F_{\Delta_i}$ we get 
\begin{equation}\label{eq:majoGammahat}
\Vvert\Gamma_\Delta-\widehat\Gamma_{\widehat\Delta}\Vvert\leq \Vvert\widehat\Gamma_\Delta-\Gamma_\Delta\Vvert+\Vvert\widehat\Gamma_{\widehat\Delta}-\widehat\Gamma_{\Delta}\Vvert.
\end{equation}
First, we deal with the first term on the right-hand side of \eqref{eq:majoGammahat}. We get classically,
\begin{equation*}
\widehat\Gamma_\Delta-\Gamma_\Delta=\frac1n\sum_{i=1}^nF_{\Delta_i}\otimes F_{\Delta_i}-\mathbb E[F_\Delta\otimes F_\Delta].
\end{equation*}
We remark that $F_{\Delta_i}\otimes F_{\Delta_i}-\mathbb E[F_\Delta\otimes F_\Delta]$, $i=1,\hdots,n$ is an i.i.d. sequence of centered random variables taking values in the space of Hilbert-Schmidt operators equipped with its usual norm $\|\cdot\|_{HS}$ defined by $\|S\|_{HS}^2=\sum_{j\geq 1}\| Se_j\|^2$ for $(e_j)_{j\geq 1}$ an orthonormal (hilbertian) basis of $\mathbb L^2([0,1])$ (recall that the definition of the norm does not depend on the basis).  Furthermore, we have
\[
\Vvert\widehat\Gamma_\Delta-\Gamma_\Delta\Vvert\leq\|\widehat\Gamma_\Delta-\Gamma_\Delta\|_{HS}
\]
 and
\[
\|\widehat\Gamma_\Delta-\Gamma_\Delta\|_{HS}^2= \frac1{n^2}\sum_{i,j=1}^n\langle F_{\Delta_i}\otimes F_{\Delta_i}-\mathbb E[F_\Delta\otimes F_\Delta],F_{\Delta_j}\otimes F_{\Delta_j}-\mathbb E[ F_\Delta\otimes F_\Delta]\rangle_{HS},
\]
where, for two Hilbert-Schmidt operators $S$ and $T$, the Hilbert-Schmidt scalar product is defined by $\langle S_1,S_2\rangle_{HS}=\sum_{j\geq 1}\langle S_1e_j,S_2e_j\rangle$ and does not depend on the considered basis. Then, taking $(e_j)_j=(\eta_j)_j$, we remark that
\begin{eqnarray*}
&&\hspace{-1.5cm}\mathbb E[\langle F_{\Delta_i}\otimes F_{\Delta_i}-\mathbb E[ F_\Delta\otimes F_\Delta], F_{\Delta_j}\otimes F_{\Delta_j}-\mathbb E[ F_\Delta\otimes F_\Delta]\rangle_{HS}]\\
&=&\sum_{\ell\geq1}\mathbb E\left[\langle( F_{\Delta_i}\otimes F_{\Delta_i}-\mathbb E[ F_\Delta\otimes F_\Delta])\eta_\ell,( F_{\Delta_j}\otimes F_{\Delta_j}-\mathbb E[ F_\Delta\otimes F_\Delta])\eta_\ell\rangle\right]\\
&=&\sum_{\ell\geq1}\int_0^1{\rm Cov}\left( F_{\Delta_i}\otimes F_{\Delta_i}\eta_\ell(t), F_{\Delta_j}\otimes F_{\Delta_j}\eta_\ell(t)\right)dt\\
&=&\mathbf1_{i=j}\sum_{\ell\geq1}\int_0^1{\rm Var}\left( F_\Delta\otimes F_{\Delta}\eta_\ell(t)\right)dt\\
&=&\mathbf1_{i=j}\sum_{\ell\geq1}\int_0^1{\rm Var}\left(\langle F_\Delta,\eta_\ell\rangle F_\Delta(t)\right)dt\\
&\leq&\mathbf1_{i=j}\sum_{\ell\geq1}\mathbb E\left[\langle F_\Delta,\eta_\ell\rangle^2\int_0^1 F_\Delta^2(t)dt\right]=\mathbf1_{i=j}\mathbb E\left[\| F_\Delta\|^4\right].
\end{eqnarray*}
Then, finally, we get
\begin{equation}\label{eq:majoGammahat1}
\mathbb E\left[\Vvert\widehat\Gamma_\Delta-\Gamma_\Delta\Vvert^2\right]\leq\frac{\mathbb E\left[\|F_\Delta\|^4\right]}{n}.
\end{equation}
Then for the second term on the right-hand side of \eqref{eq:majoGammahat},
\[
\Vvert\widehat\Gamma_\Delta-\widehat\Gamma_{\widehat\Delta}\Vvert=\Vvert(F_{\widehat W}-F_W)\otimes(F_{\widehat W}-F_W)\Vvert=\|F_{\widehat W}-F_{W}\|^2.
\]
Then, since $F_\Delta(t) = F_\Pi(t)-\mathbb E[F_\Pi(t)]$,
\begin{eqnarray}
\mathbb E[\Vvert\widehat\Gamma_\Delta-\widehat\Gamma_{\widehat\Delta}\Vvert]&=&\mathbb E[\|F_{\widehat W}-F_W\|^2]\nonumber\\
&=&\mathbb E\left[\int_0^1\left(F_{\widehat W}(t)-F_W(t)\right)^2dt\right]\nonumber\\
&=&\int_0^1{\rm Var}\left(F_{\widehat W}(t)\right)dt\nonumber\\
&=&\int_0^1{\rm Var}\left(\frac1n\sum_{i=1}^nF_{\Pi_i}(t)\right)dt\nonumber\\
&=&\frac1n\int_0^1{\rm Var}\left(F_{\Pi}(t)\right)dt=\frac1n\int_0^1\mathbb E[F_{\Delta}(t)^2]dt= \frac1n\mathbb E[\|F_\Delta\|^2].\label{eq:majoGammahat2}
\end{eqnarray}
Finally, combining eq.~\eqref{eq:majoGammahat}, \eqref{eq:majoGammahat1} and \eqref{eq:majoGammahat2}, we get
\begin{equation}\label{eq:majo_diffop}
\mathbb E[\Vvert\Gamma_{\Delta}-\widehat\Gamma_{\widehat\Delta}\Vvert^2]\leq4 \frac{\mathbb E\left[\|F_\Delta\|^4\right]}{n},
\end{equation}
and Inequality~\eqref{lambda-est} of the main file follows from Bosq inequality~\eqref{eq:Bosqeiv}.

From \cite{Bosq}, Lemma 4.3, p.104, we have the upper-bound 
\[
\|\widehat\eta_j-\widetilde\eta_j\|\leq 2\sqrt{2}\delta_j^{-1}\Vvert\Gamma_{\Delta}-\widehat\Gamma_{\widehat\Delta}\Vvert,
\]
and Inequality~\eqref{eta-est} of the main file follows from~\eqref{eq:majo_diffop}.
Now, we also have, by Cauchy-Schwarz's inequality, for all $\varphi\in\mathcal H_0^1$, 
\[
\left|\langle \widehat\mu_j-\widetilde\mu_j,\varphi\rangle\right|^2 = \left|\langle \widehat\eta_j-\widetilde\eta_j,\varphi'\rangle \right|^2\leq \|\widehat\eta_j-\widetilde\eta_j\|^2 \|\varphi'\|^2. 
\]
Then 
\[
\left\| \widehat\mu_j-\widetilde\mu_j\right\|^2_{\mathcal H^{-1}}\leq \|\widehat\eta_j-\widetilde\eta_j\|^2
\]
and Inequality~\eqref{mu-est}  of the main file follows. 

\bibliographystyle{imsart-nameyear}

\end{document}